\documentclass[a4paper,10pt,leqno,twoside]{article}

\usepackage[english]{babel}
\usepackage{umlaute, fancyheadings,  amsmath, amssymb, amsthm, latexsym, graphicx,subfig}

\newcommand{\ld}{\ensuremath{,\ldots,}}
\newcommand{\ssq}{\ensuremath{\subseteq}}
\newcommand{\smin}{\ensuremath{\setminus}}
\newcommand{\eps}{\ensuremath{\varepsilon}}

\newcommand{\T}{\ensuremath{\mathbb{T}}}

\newcommand{\N}{\ensuremath{\mathbb{N}}} 
\newcommand{\R}{\ensuremath{\mathbb{R}}}
\newcommand{\Z}{\ensuremath{\mathbb{Z}}}
\newcommand{\Q}{\ensuremath{\mathbb{Q}}}

\newcommand{\supp}{\ensuremath{\mathrm{supp}}}

\newcommand{\closure}{\ensuremath{\mathrm{cl}}}

\newcommand{\kreis}{\ensuremath{\mathbb{T}^{1}}}
\newcommand{\ntorus}[1][2]{\ensuremath{\mathbb{T}^{#1}}}

\newcommand{\alphlist}{\begin{list}{(\alph{enumi})}{\usecounter{enumi}}}
\newcommand{\romanlist}{\begin{list}{(\roman{enumi})}{\usecounter{enumi}}}
\newcommand{\listend}{\end{list}}

\newcommand{\mcup}{\ensuremath{\bigcup_{m\in\N}}}

\newcommand{\ikcup}{\ensuremath{\bigcup_{i=1}^k}}

\newcommand{\nLim}{\ensuremath{\lim_{n\rightarrow\infty}}}

\newcommand{\epslim}{\ensuremath{\lim_{\varepsilon\rightarrow 0}}}

\newcommand{\ntel}{\ensuremath{\frac{1}{n}}}

\newcommand{\halb}{\ensuremath{\frac{1}{2}}}

\newcommand{\viertel}{\ensuremath{\frac{1}{4}}}

\newcommand{\dreiviertel}{\ensuremath{\frac{3}{4}}}

\newcommand{\thx}{\ensuremath{(\theta,x)}}
\newcommand{\thom}{\ensuremath{\theta + \omega}}

\newtheorem{definition}{Definition}[section]
\newtheorem{thm}[definition]{Theorem}

\newtheorem{lem}[definition]{Lemma}
\newtheorem{cor}[definition]{Corollary}
\newtheorem{prop}[definition]{Proposition}

\newtheorem{bem}[definition]{Remark}

\numberwithin{equation}{section}
\newcommand{\rhofib}{\ensuremath{\rho_{\mathrm{fib}}}}
\newcommand{\htop}{\ensuremath{h_{\mathrm{top}}}}
\newcommand{\Emon}{\ensuremath{{\cal E}_\mathrm{mon}}}
\newcommand{\foot}{\begin{footnote}}

  \title{\Large\textsc{Strangely Dispersed Minimal Sets in the Quasiperiodically
      Forced Arnold Circle Map}} \author{P.A.~Glendinning\thanks{University of
      Manchester. Email: {\tt p.a.glendinning@manchester.ac.uk}},
    T.~J\"ager\thanks{Coll\`ege de France, Paris. Email: {\tt
        tobias.jager@college-de-france.fr}} and J.~Stark\thanks{Imperial College
      London. Email: {\tt j.stark@imperial.ac.uk}}}

\begin{document}

\setlength{\oddsidemargin}{0.02\textwidth}
\setlength{\evensidemargin}{0.02\textwidth}
%\setlength{\textheight}{0.9\textheight}
%\setlength{\topmargin}{-1cm}

%%%%%%%%%%%%%%%%%%%%%%%%%%%%%%%%%%%%%%%%%%%%%%%%%%%%%%%%%%%%%%%%
%%%%%%%%%%%%%%%%%%%%%%%%%%%%%%%%%%%%%%%%%%%%%%%%%%%%%%%%%%%%%%%%
%%%%%%%%%%%%%%%%%%%%%%%%%%%%%%%%%%%%%%%%%%%%%%%%%%%%%%%%%%%%%%%%

\maketitle \abstract{We study quasiperiodically forced circle endomorphisms,
homotopic to the identity, and show that under suitable conditions these
exhibit uncountably many minimal sets with a complicated structure, to
which we refer to as `strangely dispersed'. Along the way, we generalise
some well-known results about circle endomorphisms to the uniquely
ergodically forced case. Namely, all rotation numbers in the rotation
interval of a uniquely ergodically forced circle endomorphism are realised on
minimal sets, and if the rotation interval has non-empty interior then the
topological entropy is strictly positive. The results apply in particular to
the quasiperiodically forced Arnold circle map, which serves as a paradigm
example. }

%%%%%%%%%%%%%%%%%%%%%%%%%%%%%%%%%%%%%%%%%%%%%%%%%%%%%%%%%%%%%%%%
%%%%%%%%%%%%%%%%%%%%%%%%%%%%%%%%%%%%%%%%%%%%%%%%%%%%%%%%%%%%%%%%
%%%%%%%%%%%%%%%%%%%%%%%%%%%%%%%%%%%%%%%%%%%%%%%%%%%%%%%%%%%%%%%%

\section{Introduction}

Quasiperiodically forced circle (QPF) maps such as the QPF forced Arnold map $f
: \kreis \times \kreis \to \kreis \times \kreis$
\begin{equation} \label{e.arnold} f(\theta,\varphi) \ = \ \left(\thom,\varphi +
    \tau + \frac{\alpha}{2\pi}\sin(2\pi \varphi) + \beta\sin(2\pi\theta) \bmod
    1\right) \ ,\end{equation}where $\kreis = \R/\Z$ denotes the circle and
$\omega\notin{\bf Q}$, have been studied by a number of authors. The motivation
for this comes from two related directions. First, Grebogi \textit{et al}
\cite{GOPY} showed that it is possible to have strange (i.e. geometrically
complicated) nonchaotic attractors (SNAs) over a range of parameter values with
positive measure, and later (e.g. \cite{DGO,RBOAG,RBOAG2,RO}) that maps such as
\ref{e.arnold} are good candidates for simple invertible examples of such
behaviour. This aspect has been followed up in the work of Feudel and Pikovsky
and ourselves amongst others \cite{FKP,PF,Chastell,GFPS}. Secondly, from a
different perspective Herman \cite{herman:1983} had already proved the existence
of SNA in certain parameter families of QPF circle diffeomorphisms that are
induced by the projective action of $SL(2,\R)$-cocycles over an irrational
rotation.

Despite considerable interest over subsequent years, rigourous mathematical
results remained rare. The original Grebogi \textit{et al} example \cite{GOPY}
was a non-invertible map with a special structure. This structure was abstracted
by Keller, who proved the existence of SNAs under simple conditions in this
class of maps \cite{Keller}.  J{\"a}ger \cite{jaeger:2007} further analysed the
structure of such invariant sets. Subsequently, Glendinning \textit{et al}
\cite{GJK} proved that although non-chaotic in the sense of Lyapunov exponents,
such systems did exhibit sensitive dependence on initial conditions. Meanwhile,
Stark \cite{Berlin-proc} showed that SNAs in QPF maps could not be non-smooth
graphs, but had to have a more complex structure, and an extension of this
approach by Sturman and Stark \cite{Semi-uniform} showed that the normal
Lyapunov exponents of a SNA could not all be negative. Finally, new methods were
established quite recently by Bjerkl\"ov and J\"ager, which allow to prove the
existence of SNA in much broader classes of quasiperiodically forced maps than
the two mentioned above \cite{bjerkloev:2005a,jaeger:2006b,jaeger:2006c}.

Additional properties of invertible circle maps were derived in
\cite{SFGP,Monotone}, and used together with results in
\cite{jaeger/keller:2006} to give a generalization of the Poincare
classification of circle homeomorphisms \cite{jaeger/stark:2006}. Further,
J{\"a}ger and Keller \cite {jaeger/keller:2006} showed that if a QPF circle
homeomorphism, homotopic to the identity, with appropriate conditions on its
rotation number, was topologically transitive then any minimal set was
`strangely dispersed' (see below for definition). Dynamics of this type are
constructed in \cite{beguin/crovisier/jaeger/leroux:2006b}. However, it is also
known that the minimal set is unique in this situation (\cite{huang/yi:2007} or
\cite{beguin/crovisier/jaeger/leroux:2006b}), such that there is no co-existence
as in Theorem~\ref{t.sdsm}. Further, the examples given in
\cite{beguin/crovisier/jaeger/leroux:2006b} only have low regularity, and it is
still completely open whether the same phenomenon can occur in smooth systems as
well (e.g.\ in QPF analytic circle diffeomorphisms).

Here we turn to examine the behaviour of QPF maps such as the forced Arnold map
(\ref{e.arnold}) above. This is motivated both by the considerable volume of
numerical work, and the fact that the unforced Arnold map has a rich and
interesting structure has been described in some detail by MacKay and Tresser
\cite{MacKay}. This gave a beautiful description of the transition to chaotic
behaviour in the unforced case. Numerical experiments have suggested that in the
QPF map the appearance of strange nonchaotic structures occurs at the complex
boundary between the regular and chaotic parameter regions.

Unfortunately, MacKay and Tresser's analysis made heavy use of periodic orbits
and doubling cascades. Since (\ref{e.arnold})  has no periodic orbits (this follows
immediately from the fact that $\omega$ is irrational) it is not immediately
clear how to generalize their work to the QPF case. Indeed, almost all of our
understanding of chaos is based on generalizations of the horseshoe (e.g.
\cite{katok/hasselblatt:1997}) and horseshoes imply the existence of periodic
orbits, so either horseshoes are irrelevant for the study of chaos in
quasiperiodically forced systems or the chaos is essentially a suspension of a
horseshoe. If the former is the case then it is natural to ask which orbits form
the backbone of the chaos, i.e. which orbits play the role of the periodic
orbits in the horseshoe?  There are therefore at least three reasons for
considering noninvertible quasiperiodically forced circle maps, $k >1$ in
(\ref{e.arnold}). First in an attempt to obtain some rigorous results on complex
invariant sets, second as an extension of the results for noninvertible circle
maps, and third as a move towards understanding chaotic sets which are not
modelled by horseshoes. Our main motivation has been the first of these. We
shall prove that if $k$ and $\beta$ are sufficiently large then the the forced
Arnold map (\ref{e.arnold}) exhibits uncountably many minimal sets with a
complicated structure, to which we refer to as `strangely dispersed'.  In the
proof of this result it becomes necessary to prove analogues of a number of
results for noninvertible circle maps in the context of noninvertible
quasiperiodically forced circle maps. An appealing, albeit heuristic,
interpretation of this result is that in moving along a path in parameter space
from a nonchaotic state of an invertible quasiperiodically forced circle map to
a chaotic noninvertible circle map of the type discussed below it is necessary
to create strange nonchaotic invariant sets. One way of achieving this is to
create this set as a stable set, which later loses stability. If this is the
case then it goes some way towards explaining why nonchaotic strange attractors
must exist in such families.

\noindent

{\bf Acknowledgements.} We would like to thank Sylvain Crovisier for pointing
out to us the result by Bowen \cite{bowen:1971} and its consequences for the
entropy of QPF monotone circle maps. Tobias J\"ager was supported by a
research fellowship of the German Research Foundation (DFG, Ja 1721/1-1).

%%%%%%%%%%%%%%%%%%%%%%%%%%%%%%%%%%%%%%%%%%%%%%%%%%%%%%%%%%%%%%%%
%%%%%%%%%%%%%%%%%%%%%%%%%%%%%%%%%%%%%%%%%%%%%%%%%%%%%%%%%%%%%%%%
%%%%%%%%%%%%%%%%%%%%%%%%%%%%%%%%%%%%%%%%%%%%%%%%%%%%%%%%%%%%%%%%

\section{Main Results}

Let $\kreis = \R/\Z$ denote the circle and suppose $\Theta$ is a compact
metric space and $r:\Theta \to \Theta$ a continuous map. We consider
skew-products on $\Theta \times \kreis$ given by continuous maps $f : \Theta
\times \kreis \to \Theta \times \kreis$ of the form
\begin{equation} \label{e.f} f(\theta,\varphi) \ = \ (r(\theta),f_\theta(\varphi)) \ .
\end{equation}
The case we are primarily interested in is that of QPF
circle endomorphisms, that is $\Theta = \kreis$ and $r(\theta) =
\theta+\omega$ with $\omega \in \kreis$ irrational. However, some of the
results we obtain naturally generalise to the uniquely ergodically forced
(UEF) case (meaning that there exists a unique
$r$-invariant probability measure $\mu$ on $\Theta$).

The maps $f_\theta$ in (\ref{e.f}) will be called {\em fibre maps}. Most of
the time, we will assume in addition that $f$ is homotopic to the map $(\theta,\varphi)
\mapsto (r(\theta),\varphi)$. If this holds, we say $f$ is {\em homotopic to the
identity} (slightly abusing terminology in the case that $r$ is not homotopic to the
identity on $\Theta$). Then all $f_\theta$ are circle endomorphisms of degree
one, and further there exist a continuous lift $F:\Theta \times \R \to \Theta
\times \R$ that satisfies $\pi \circ F = f \circ \pi$, where $\pi : \Theta
\times \R \to \Theta \times \kreis$ denotes the natural projection. Moreover,
if $\Theta$ is connected, then these lifts are always uniquely defined modulo
an integer. In the same way we can define the continuous
lifts $F_\theta:\R\to \R$ of fibre maps $f_\theta$ which satisfy
$\pi \circ F_\theta = f_\theta \circ \pi$, where $\pi : \R \to \kreis$
denotes the natural projection, with the obvious
abuse of notation on the projection operators $\pi$.

We define the {\em fibred rotation interval} of
a lift $F$ by
\begin{equation} \label{e.rhofib} \rhofib(F) \ := \ \left\{
    \left. \limsup_{n\to\infty} \frac{1}{n}(F^n_\theta(x) -x)\right| \thx \in
    \Theta \times \R \right\} \ .
\end{equation}
where $F^n_\theta (x)=F_{r^{n-1}(\theta )} \circ \ldots \circ
F_\theta(x)$. Observe that $\rhofib$ for two different lifts of the same UEF
endomorphism of $ \Theta \times \kreis$ will differ by an integer translation.

An important special case will be the one of UEF (or QPF) {\em monotone circle
  maps}, by which we mean that each fibre map $f_\theta$ preserves the cyclic
order on $\kreis$ (but we allow $f$ to be non-injective). This is true if and
only if the fibre maps $F_\theta$ of any lift $F$ of $f$ are monotonically
increasing. It is a well-known result of Herman \cite{herman:1983} that the
fibred rotation interval of a UEF monotone circle map is always a single point
(restated below as Theorem~\ref{t.herman}).

\begin{bem}
  Note that in the QPF case there are in general several ways of assigning a
  rotation set to a torus endomorphism which is homotopic to the identity, as
  discussed very concisely in \cite{misiurewicz/ziemian:1989}. However, if $f$
  has skew-product structure as in (\ref{e.f}), then all these different notions
  coincide. This follows easily from Theorem~\ref{t.rotation-interval} below, in
  combination with \cite[Theorem 2.4 and Corollary
  2.6]{misiurewicz/ziemian:1989}. The above definition is the one which is most
  convenient for our purposes, and we have adapted it to the fibred setting by
  projecting to the second coordinate, thus obtaining a subset of the real line
  instead of a subset of $\R^2$ for a general endomorphism of \ntorus[2].
\end{bem}

Recall that a closed, $f$-invariant set $M$ is {\em minimal} if it contains no
proper $f$-invariant closed subset \cite{katok/hasselblatt:1997}. This is
equivalent to the orbit of every point in $M$ being dense in $M$. The {\em
  topological entropy} $\htop(f)$ of a map $f$ is a common measure of the
complexity of its dynamics, and indeed provides one of the standard definitions
of chaotic behaviour \cite{katok/hasselblatt:1997}. A definition and brief
overview is given below in Section~\ref{Entropy}. The following theorem is then
a generalisation of well-known results on unforced circle endomorphisms (see,
for example, \cite{Misnew}).

\begin{thm} \label{t.rotation-interval} Suppose $F$ is the lift of a UEF circle
  endomorphism $f : \Theta \times \kreis \to \Theta \times \kreis$, homotopic to
  the identity. Then $\rhofib(F)$ is a closed interval (including the
  possibility of a singleton $\rhofib(F) = [\rho,\rho]$). For any $\rho \in
  \rhofib(F)$ there exists a minimal set $M_\rho \subset \Theta \times \kreis$
  with the following properties: \romanlist
\item $\frac{1}{n}(F^n_\theta(x) -x)$ converges uniformly to $\rho$ on 
$\pi^{-1}(M_\rho)$ as $n \to \infty$.
\item $\htop(f|_{M_\rho}) = 0$.
\listend
\end{thm}

The proof of (i) is given in Section~\ref{PlateauMaps} and that of (ii) at the end
of Section~\ref{Entropy}. Although the dynamics on each $M_\rho$ is simple, if
the rotation interval is non-trivial, the overall dynamics of the map is
complex:

\begin{thm} \label{t.entropy} Suppose $f$ is a UEF circle endomorphism,
  homotopic to the identity, with lift $F$. If $\rhofib(F)$ has non-empty
  interior, then $\htop(f) > 0$.
\end{thm}

The proof is given in Section~\ref{Entropy}. We remark that for a QPF monotone
circle map $f$ the situation is quite different. As mentioned above, the
rotation interval is reduced to a single point in this case, and the topological
entropy is always zero. The latter is a more or less direct consequence of a
result by Bowen \cite{bowen:1971}, see Section~\ref{Entropy} below.
\medskip

Once these basic facts are established, we can turn to a new phenomenon which is
specific to the quasiperiodically forced setting. In the case of unforced circle
endomorphisms, minimal sets may be either periodic orbits or Cantor sets,
corresponding to rational and irrational rotation numbers, respectively. In the
quasiperiodically forced case however, they can have a much more complicated
structure. In order to make this precise, we introduce the following notion.

\begin{definition} \label{d.sdms} Suppose $f$ is a QPF circle endomorphism,
  homotopic to the identity. We say a compact subset $M\ssq \ntorus$ is a
  strangely dispersed minimal set, if it has the following three properties:
  \romanlist
\item $M$ is a minimal set.
\item Every connected component $C$ of $M$ is contained in a single fibre, that
  is $\pi_1(C)$ is a singleton.
\item For any point $\thx \in M$ and any open neighbourhood $U$ of $\thx$, the
  set $\pi_1(U\cap M)$ contains a non-empty open interval. \listend
\end{definition}

\begin{bem}~ \alphlist
\item Property (iii) is a rather direct consequence of (i) (see
  Section~\ref{SDMS}). We have only included it here to emphasize the
  peculiarity of property (ii).
\item It is actually not difficult to construct a set which has properties (ii)
  and (iii). Indeed, let $(a_\eta)_{\eta\in \kreis \cap \Q}$ be any sequence of
  strictly positive real numbers with $\sum_{\eta\in\kreis \cap \Q} a_\eta =
  1$. For any $\theta \in \kreis$, let $\phi(\theta) := \sum_{\eta \in
    [0,\theta]\cap \Q} a_\eta$. Then the topological closure of the graph $\Phi
  := \{(\theta,\phi(\theta)) \mid \theta \in \kreis\}$ of $\phi$ is a compact
  set that has these two properties. Of course, the interesting point in the
  above definition is to have a set with this structure as the minimal set of a
  dynamical system. \listend
\end{bem}

It will follow from our arguments in Section~\ref{SDMS} that the appearance of
strangely dispersed minimal sets is a rather general phenomenon for QPF circle
endomorphisms, provided that the quasiperiodic forcing has a certain
strength. However, for simplicity we will formulate the results only for a
particular example, namely for the QPF Arnold Circle Map (\ref{e.arnold}).

\begin{thm} \label{t.sdsm} Suppose $f$ is given by (\ref{e.arnold}), with
  driving frequency $\omega \in \kreis\smin\Q$ and real parameters $\tau,
  \alpha$ and $\beta$.

\alphlist
\item If $\alpha > 1$ and $|\beta| \geq \frac{3}{2}$, then for any
  $\rho\in\rhofib(f)$ there exists a strangely dispersed minimal set $M_\rho$
  which satisfies properties (i) and (ii) in Theorem~\ref{t.rotation-interval}~.
\item If $|\alpha| \geq \frac{5}{2}\pi$, then $\rhofib(f)$ has length $\geq 1$,
  in particular its interior is non-empty. \listend
\end{thm}

\begin{bem}~
\alphlist

\item The bounds given here are not optimal and may surely be improved by a more
careful analysis. Further, part (b) of this theorem is rather trivial, but it is
important for the interpretation of (a). Namely, if both conditions in
(a) and (b) are satisfied simultaneously, we obtain the existence of
uncountably many pairwise disjoint and strangely dispersed minimal sets, one
for each rotation number in the rotation interval. Albeit most likely
superficial, this presents an intriguing analogy to the theory of twist maps,
where at suitable parameter values the standard map exhibits uncountably many
Aubry-Mather sets, again one for each rotation number in the rotation
interval.

\item As indicated in the introduction above, the existence of strangely
  dispersed minimal sets is already known in QPF circle homeomorphisms
  \cite{jaeger/keller:2006, beguin/crovisier/jaeger/leroux:2006b, huang/yi:2007}
  though existing constructions only work in maps of low regularity.

\listend
\end{bem}

%%%%%%%%%%%%%%%%%%%%%%%%%%%%%%%%%%%%%%%%%%%%%%%%%%%%%%%%%%%%%%%%
%%%%%%%%%%%%%%%%%%%%%%%%%%%%%%%%%%%%%%%%%%%%%%%%%%%%%%%%%%%%%%%%
%%%%%%%%%%%%%%%%%%%%%%%%%%%%%%%%%%%%%%%%%%%%%%%%%%%%%%%%%%%%%%%%

\section{Plateau Maps and Rotation Numbers: Proof of Theorem~\ref{t.rotation-interval}}
\label{PlateauMaps}

In order to prove Theorem~\ref{t.rotation-interval}, we will first be
concerned with unforced circle endomorphisms and their lifts. The basic idea,
which is the use of plateau maps to identify orbits with a given rotation
number, was first introduced by Boyland \cite{boyland:} (see also
\cite{Misnew} for a survey). Let ${\cal E}$ denote the space of continuous
maps $G:\R \to \R$ which are the lift of a circle endomorphism of degree
one. The latter just amounts to saying that
\begin{equation}\label{e.periodicity}
G(x+k) \ = \ G(x)+k \quad \forall k\in\Z \ .
\end{equation}
We equip ${\cal E}$ with the topology of uniform convergence. Further, we
denote by $\Emon$ the space of all maps in $G \in {\cal E}$ which are
monotonically increasing. Then $G \in \Emon$ if and only if it is the lift of
a monotone circle map of degree one. Note that we explicitly allow for $G \in
\Emon$ to be non-injective, in which case there exist intervals that are
mapped to a single point by $G$. We refer to such intervals as {\em
plateaus}, and call maps in \Emon\ {\em plateau maps} (including the case
when there are no plateaus, for simplicity).

\begin{figure}[h]
\centering
\subfloat{
\includegraphics{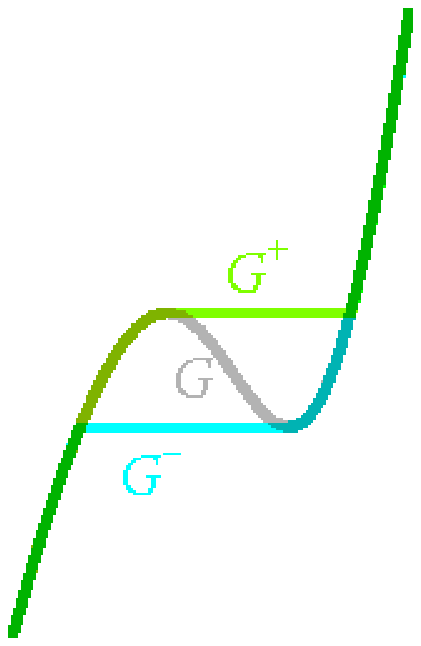}
}
\quad
\subfloat{
\includegraphics{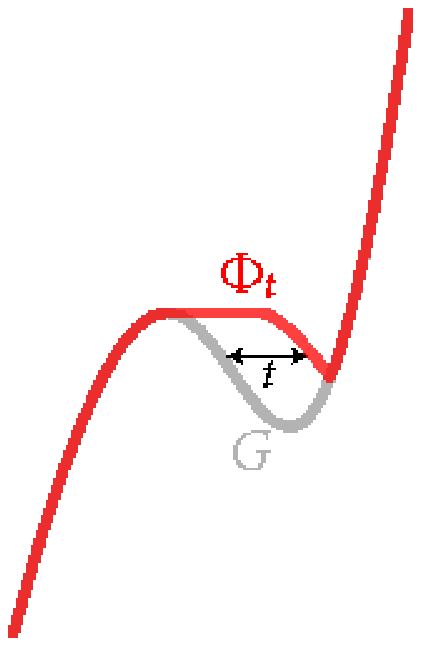}
}
\quad
\subfloat{
\includegraphics{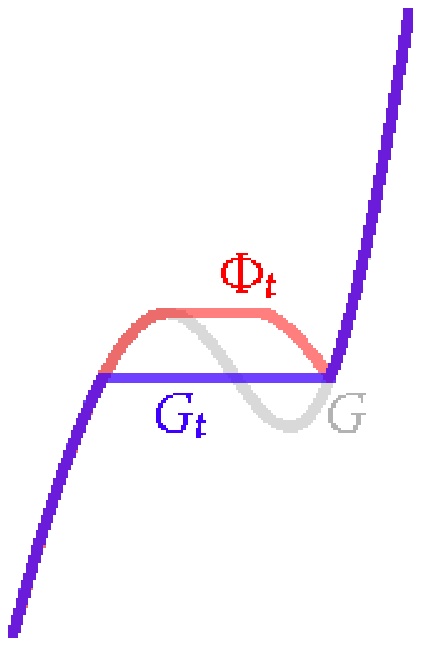}
}
\caption{Illustration of  the plateau maps $G^+$ and $G^-$ and the functions $\Phi_t$ and $G_t$. }
\label{Phi_fig}
\end{figure}

For any $G\in\Emon$, let ${\cal U}(G)$ denote the union of the interiors of all
plateaus of $G$. In other words
$$
{\cal U}(G) \ = \ \{ x \in \R \mid \exists \eps > 0 : G([x-\eps,x+\eps]) = \{G(x) \} \} \ .
$$
Now suppose $G \in {\cal E}$. We assign to $G$ a pair of plateau maps $G^-
\leq G \leq G^+$ (Figure \ref{Phi_fig}) by
$$
G^+(x) \ := \ \sup_{\xi \leq x} G(\xi) \quad \textrm{and} \quad G^-(x) \ := \ \inf_{\xi \geq x} G(\xi) \ .
$$

Note that if $G$ is a plateau map itself, then $G^- = G^+ = G$. Further, it
follows easily from (\ref{e.periodicity}) that 
\begin{equation}\label{finite_interval}
G^+(x) = \sup_{\xi \in [x-1,x]} \quad \textrm{and} \quad G^-(x) = \inf_{\xi \in [x,x+1]} G(\xi) \ .
\end{equation}
The reason why plateau maps are such a convenient tool for the computation of
rotation intervals is the fact that they always have a uniquely defined
rotation number, and this remains true in the quasiperiodically forced setting
(see Theorem~\ref{t.herman} below). Furthermore, as the following proposition
shows, there always exists a homotopy between the maps $G^-$ and $G^+$ with
some additional nice properties, and this will be the key ingredient in the
proof of Theorem~\ref{t.rotation-interval}~.

\begin{prop} \label{p.plateau-homotopy} There exists a continuous mapping $\R
  \times [0,1] \times {\cal E} \to \R$, $(x,t,G) \mapsto G_t(x)$, with the
  following properties: \romanlist
\item The family $(G_t)_{t\in[0,1]}$ is a homotopy between $G^-$ and $G^+$, that
  is $G_0 = G^-$ and $G_1 = G^+$.
\item For all $t\in[0,1]$ we have $G_t \in \Emon$.
\item For all $x \in \R$ the map $t\mapsto G_t(x)$ is
monotonically increasing.
\item If $G_t(x) \neq G(x)$, then $x \in {\cal U}(G_t)$.
\listend
\end{prop}
Note that due to the periodicity property (\ref{e.periodicity}) of $G\in
{\cal E}$ and compactness, the induced mapping $[0,1]\times {\cal E} \to \Emon,\
(t,G) \mapsto G_t$ is continuous as well.
\proof\ The mappings $\R \times {\cal E} \to \R$, $(x,G) \mapsto G^\pm(x)$ are
clearly continuous and monotone in $x$. They are also continuous and
monotonically increasing in $G$, the latter with respect to the partial ordering
on ${\cal E}$ given by $G_1 \leq G_2$ if $G_1(x) \leq G_2(x)\ \forall
x\in\R$. Similarly, the mapping (Figure \ref{Phi_fig})
\begin{equation} \label{e.phit}
\Phi : \R \times [0,1] \times {\cal E} \to \R \quad , \quad (x,t,G) \mapsto
\Phi_t(x) := \sup_{\xi \in [x-t,x]} G(\xi) \ ,
\end{equation}
is continuous and monotonically increasing $t$ and $G$. We define our required
homotopy by
\begin{equation} \label{e.Gt}
G_t(x) \ := \ (\Phi_t)^-(x) \ = \ \inf_{\zeta \geq x} \;  \sup_{\xi \in[\zeta-t,\zeta]} G(\xi) \ .
\end{equation}

For any given $x \in {\cal E}$ the function $(x,t) \mapsto G_t(x)$ is continuous
as the composition of continuous functions. By definition $\Phi_0 = G$, and
hence $G_0 = G^-$. Also, by (\ref{finite_interval}) we have $\Phi_1 = G^+$ and
hence $G_1 = (G^+)^- = G^+$. This proves (i). For any $t\in[0,1]$ the map $G_t =
(\Phi_t)^-$ is a plateau map, since it is in the image of the mapping ${\cal E}
\to \Emon,\ G \mapsto G^-$. Thus (ii) holds.  The monotonicity of the mapping
$\Phi$ in $t$ and of the mapping $G \mapsto G^-$ in $G$ immediately implies
(iii).

It remains to prove (iv). We first show that if for a given $t\in[0,1]$ and
$x\in \R$ we have $G_t(x) < \Phi_t(x)$ then $G_t$ is constant in an open
neighbourhood of $x$. Since $G_t = (\Phi_t)^-$, and
$\Phi_t$ is continuous, then if $G_t(x) < \Phi_t(x)$ there must exist some
$\xi_0 > x$ with $\Phi_t(\xi_0) = G_t(x)$. By the continuity of $\Phi_t$ we
further have $\Phi_t(x') > G_t(x)$ for all $x'$ in a small open neighbourhood
$U$ of $x$. Without loss of generality, we shrink $U$ so that it does not
contain $\xi_0$, which implies that $x' < \xi_0$ for all $x' \in U$. This means
that $\inf_{\xi \geq x'} \Phi_t(\xi) \leq \Phi_t(\xi_0)$ for all $x' \in U$. If
$x' \geq x$, then automatically $\inf_{\xi \geq x'} \Phi_t(\xi) \ge \inf_{\xi
  \geq x} \Phi_t(\xi) = \Phi_t(\xi_0)$, whereas if $x' < x$, then $\inf_{\xi
  \geq x'} \Phi_t(\xi) =\min \{\inf_{x' \le \xi < x} \Phi_t(\xi),\inf_{\xi \geq
  x} \Phi_t(\xi) \} = \min \{\inf_{x' \le \xi < x} \Phi_t(\xi),\Phi_t(\xi_0) \}
= \Phi_t(\xi_0)$, since $\Phi_t(\xi) > G_t(x) = \Phi_t(\xi_0)$ for all $\xi \in
U$. Hence $G_t(x') = \Phi_t(\xi_0)$ for all $x'\in U$, so that $G_t$ is constant
on $U$ as required.

\begin{figure}
\centering
\subfloat[The case $G_t(x) < G(x)$]{
\includegraphics{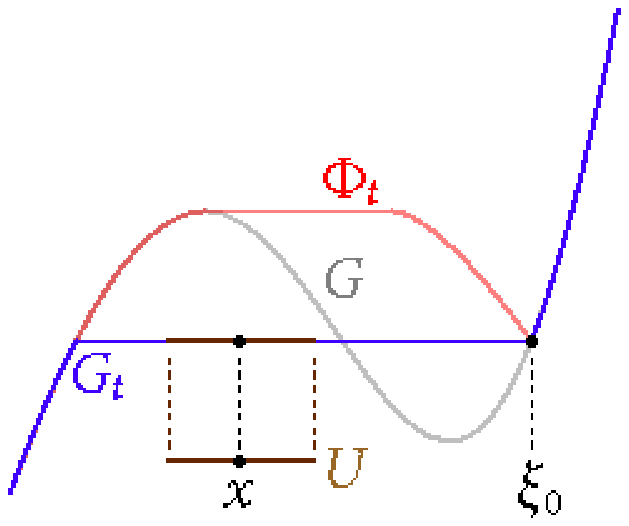}
}
\qquad
\subfloat[The case $G_t(x) > G(x)$]{
\includegraphics{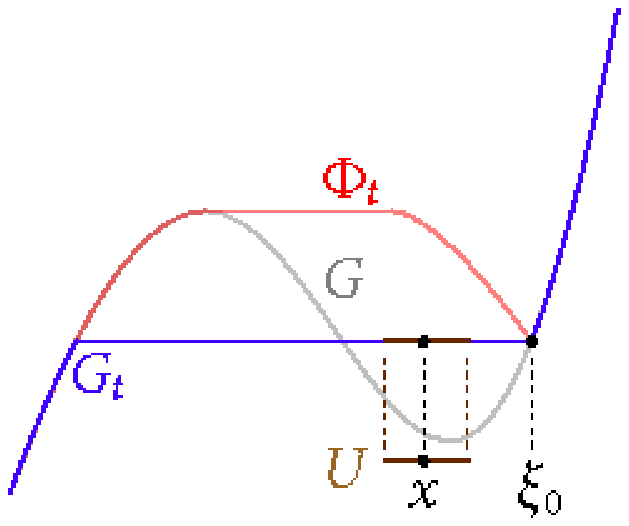}
}
\caption{Proof of Proposition \ref{p.plateau-homotopy} (iv) when $G_t(x) < \Phi_t(x)$. }
\label{Plateau_Proof_Fig}
\end{figure}

To now prove (iv), fix $t\in[0,1]$ and $x\in \R$ with $G_t(x) \neq G(x)$. First
suppose $G_t(x) < G(x)$ (Figure \ref{Plateau_Proof_Fig}a). Since $G(x) \leq
\Phi_t(x)$ for all $x\in \R$, this implies that $G_t(x) < \Phi_t(x)$ and hence
by the above $G_t$ is constant on an open neighbourhood of $x$, as required.

On the other hand, if $G_t(x) > G(x)$, we consider two cases. By definition, we
always have $G_t(x) \le \Phi_t(x)$. Hence either $G_t(x) < \Phi_t(x)$ (Figure
\ref{Plateau_Proof_Fig}b) or $G_t(x) = \Phi_t(x)$. In the former case we again
apply the argument above.

\begin{figure}
\centering
\includegraphics{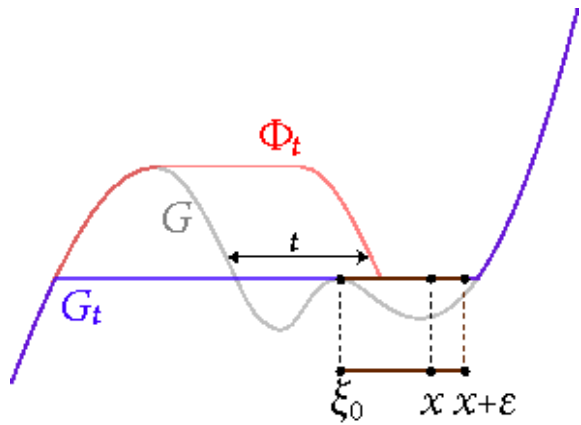}
\caption{Proof of Proposition \ref{p.plateau-homotopy} iv) for the case  $G(x) < G_t(x) = \Phi_t(x)$. }
\label{Plateau_Proof_Fig_Eq}
\end{figure}

In the latter case, $G_t(x) = \Phi_t(x)$, the definition of $\Phi_t$ implies
that there exists some $\xi_0 \in [x-t,x)$ such that $\Phi_t(x) = G(\xi_0)$
(Figure \ref{Plateau_Proof_Fig_Eq}). For any $x' \in [\xi_0,x]$ we have $\xi_0
\in [x'-t,x']$ and hence $\Phi_t(x') \geq G(\xi_0) = \Phi_t(x)$. Thus $G_t(x') =
\inf_{\xi \geq x'} \Phi_t(\xi) = \inf_{\xi \geq x} \Phi_t(\xi) = G_t(x)$, so
that $G_t$ is constant on a left neighbourhood of $x$. Now, by definition
$\Phi_t(x') \ge G_t(x)$ for all $x' \ge x$, and since $G_t(x) = \Phi_t(x)$, we
have $\Phi_t(x') \ge \Phi_t(x)$ for all $x' \ge x$. Since $G(x) < G_t(x) =
\Phi_t(x)$, the continuity of $G$ implies that there exists an $\epsilon > 0$
such that $G(x') < \Phi_t(x)$ for all $x' \in [x,x+\epsilon)$ Furthermore, by
the defintion of $\Phi_t$ we have $G(x') \le \Phi_t(x)$ for all $x' \in
[x-t,x]$. Hence for any $x' \in [x,x+\epsilon)$ we have $G(\xi) \le \Phi_t(x)$
for all $\xi \in [x'-t,x']$. Thus $\Phi_t(x') \le \Phi_t(x)$ for all $x' \in
[x,x+\epsilon)$ and so $\Phi_t$ is constant on $[x,x+\epsilon)$. By definition
$G_t$ is non-decreasing, so $G_t(x') \ge G_t(x)$ for all $x' \in
[x,x+\epsilon)$. But $G_t(x') \leq \Phi_t(x')$ for any $x'$, and so in particular
for $x' \in [x,x+\epsilon)$ we have $G_t(x') \leq \Phi_t(x') = \Phi_t(x) =
G_t(x)$. Hence $G_t(x') = G_t(x)$ for all $x' \in [x,x+\epsilon)$, and so $G_t$
is also constant on a right neighbourhood of $x$. This completes the proof of
(iv).

\qed
\medskip

The following lemma will provide the link between the orbits of the original map
and the plateau maps derived from it.
\begin{lem} \label{l.orbitlink} Suppose $(G_n)_{n\in\N_0}$ is a sequence of
  plateau maps and let $G^{(n)} := G_{n-1} \circ \ldots \circ G_0$. Then there
  exists $x\in\R$ with the property that $G^{(n)}(x) \notin {\cal U}(G_n) \
  \forall n \in \N_0$.
\end{lem}
\proof\ We argue by contradiction. Suppose for all $x\in\R$, there exists some
$n \in \N_0$, such that $G^{(n)}(x) \in {\cal U}(G_n)$. Let $V_n :=
(G^{(n)})^{-1}({\cal U}(G_n))$. Then the open sets $\pi(V_n)$ form an open
cover of \kreis\ and hence, since \kreis\ is compact, there is a finite
subcover. Thus $\kreis \ssq \pi(V_0) \cup \ldots \cup \pi(V_N)$ for some
$N\in\N_0$ and hence $\R\ssq V_0 \cup \ldots \cup V_N$. However, as every
plateau is mapped to a single point and there are at most countably many
plateaus, this implies that $G^{(N)}(\R)$ is countable and therefore a strict
subset of $\R$. Since all $G_n$ are surjective, this yields the required
contradiction.

\qed
\medskip

Now we can turn to the forced setting. Recall that we say $f$ is a UEF
monotone circle map, if all of its fibre maps $f_\theta$ are circle maps of
degree one which preserve the cyclic order on \kreis. This is true if and only
if any lift $F:\Theta \times \R \to \Theta \times \R$ of $f$ satisfies
$F_\theta \in \Emon \ \forall \theta \in \Theta$. As mentioned before, the
rotation number of UEF monotone circle maps is uniquely defined.
\begin{thm}[Herman \cite{herman:1983, SFGP}] \label{t.herman} Suppose $f$ is a
  UEF monotone circle map, homotopic to the identity, with lift $F$. Then the
  limit
\begin{equation} \label{e.rotnum}
\rho(F) \ := \ \nLim \frac{1}{n}(F^n_\theta(x) - x)
\end{equation}
exists and is independent of $\thx$, and the convergence in (\ref{e.rotnum}) is
uniform on $\Theta \times \R$. Furthermore, $\rho(F)$ depends continuously on
$F$. We call $\rho(F)$ the {\em fibred rotation number} of $F$.
\end{thm}
In fact, the result in \cite{herman:1983} is only stated for UEF circle
homeomorphisms, but the proof given there literally goes through in this
slightly more general situation. Alternatively, \cite{SFGP} explicitly proves
the existence of a unique rotation number for non-strictly monotone maps.

\proof[\bf Proof of Theorem~\ref{t.rotation-interval} (i) ] Suppose $f : \Theta
\times \kreis \to \Theta \times \kreis, ( \theta,\varphi) \mapsto
(r(\theta),f_\theta(\varphi))$ is a UEF circle endomorphism homotopic to the
identity and $F \in {\cal E}$ is a lift of $f$. We define two UEF monotone maps
$F^-,F^+$ by $F^-_\theta := (F_\theta)^-$ and $F^+_\theta := (F_\theta)^+$. Then
$F^-$ and $F^+$ are the lifts of two UEF monotone circle maps, and by
Theorem~\ref{t.herman} the fibred rotation numbers of $F^-$ and $F^+$ are
well-defined. Since $F^-_\theta(x) \leq F_\theta(x) \leq F^+_\theta(x) \ \forall
\thx \in \Theta \times \R$ it follows easily that $\rhofib(F) \ssq
[\rho(F^-),\rho(F^+)]$.

We obtain a homotopy $F_t$ from $F^-$ to $F^+$ by defining $F_{t,\theta}(x) :=
(F_\theta)_t(x)$, where $(x,t,G) \mapsto G_t(x)$ is the mapping provided by
Proposition~\ref{p.plateau-homotopy}~. Note that each $F_t$ is continuous,
because $F_\theta$ depends continuously on $\theta$ and the mapping $(x,t,G)
\mapsto G_t(x)$ is continuous. Since $t \mapsto F_t$ is continuous and
monotone (by property (iii) of the proposition), and as the fibred rotation
number depends continuously on the system, the mapping $t\mapsto \rho(F_t)$ is
a continuous and monotonically increasing function. Therefore, it maps the
interval $[0,1]$ surjectively onto $[\rho(F^-),\rho(F^+)]$. Consequently, for
any fixed $\rho\in[\rho(F^-),\rho(F^+)]$ there exists some $t = t(\rho) \in
[0,1]$, such that $\rho(F_t) = \rho$. Fixing any $\theta_0\in\Theta$ and
applying Lemma~\ref{l.orbitlink} with $G_n = F_{t,r^n(\theta_0)}$, we
obtain an $x_0\in\R$ with $F^n_{t,\theta_0}(x_0) \notin
{\cal U}(F_{t,r^n(\theta_0)})\ \forall n\in\N_0$. By property (iv)
in Proposition~\ref{p.plateau-homotopy}, we have $\{ x \in \R \mid
F_{t,\theta}(x) \neq F_\theta(x) \} \ssq {\cal U}(F_{t,\theta}) \ \forall
\theta \in \Theta$. Therefore $F^n_{t,\theta_0}(x_0) = F^n_{\theta_0}(x_0) \
\forall n\in\N_0$, which means that the orbits of $(\theta_0,x_0)$ under the
maps $F_t$ and $F$ coincide. Hence
\[
\nLim \frac{1}{n}(F^n_{\theta_0}(x_0)-x_0) \ = \ \nLim \frac{1}{n}(F^n_{t,\theta_0}(x_0)-x_0) \
= \ \rho(F_t) \ = \rho \ .
\]
This shows that $\rho$ is contained in $\rhofib(F)$, and since $\rho \in
[\rho(F^-),\rho(F^+)]$ was arbitrary we obtain $\rhofib(F) =
[\rho(F^-),\rho(F^+)]$.

Furthermore, by continuity it follows that $F_\theta(x) = F_{t,\theta}(x)$
for all $\thx$ in the set
\[
A \ := \ \closure\left(\left\{ F^n_{\theta_0}(x_0+k) \mid n \in
\N_0, k\in \Z\right\}\right) \ ,
\]
where $\closure(\cdot)$ denotes the topological closure. If we define
$M_\rho$ as the omega limit set of $\pi(\theta_0,x_0)$, that is
\[
M_\rho \ = \ \cap_{n\geq 0}
\closure(\{f^k\circ \pi(\theta_0,x_0) \mid k \geq n\}) \ ,
\]
then clearly $\pi^{-1}(M_\rho) \ssq A$. Hence the restrictions of $F$ and
$F_t$ to $\pi^{-1}(M_\rho)$ coincide. It follows that the quantities
$\frac{1}{n}(F^n_\theta(x)-x)$ converge uniformly to $\rho$ on $\pi^{-1}(M_\rho)$ as
$n$ tends to infinity, since this is true for the quantities
$\frac{1}{n}(F^n_{t,\theta}(x)-x)$ by Theorem~\ref{t.herman}~. \qed

\begin{bem}
\label{r.projections}
Note that the minimal sets $M_\rho$ of Theorem~\ref{t.rotation-interval} have to
project down to a minimal set for the underlying transformation $r:\Theta \to
\Theta$. Since we assume $r$ to be uniquely ergodic, the only such minimal set
is the topological support $\supp(\mu)$ of the unique $r$-invariant
probability measure $\mu$. Thus for any $\rho_1,\rho_2 \in \rhofib(F)$ we have $\pi_1(M_{\rho_1})
\cap \pi_1(M_{\rho_2}) = \supp(\mu) \neq \emptyset$. In particular, if
$r$ is an irrational rotation, $\pi_1(M_{\rho}) =\kreis$ for any $\rho \in \rhofib(F)$.
\end{bem}

\begin{bem}
\label{r.avoid-plateaus}
Suppose $\rho \in \rhofib(F)$ and $t=t(\rho)$ and $F_t$ are chosen as in the
above proof. Denote by $f_t$ the UEF monotone circle map induced by $F_t$. Then
the minimal sets $M_\rho$ defined in the above proof have the property that
they do not intersect the set of plateaus of $f_t$, that is $M_\rho$ is
disjoint from $\pi({\cal V}(F_t))$ where
\[
{\cal V}(F_t) := \bigcup_{\theta \in \Theta} \{\theta\} \times {\cal
U}(F_{t,\theta}) \ .
\]
By invariance, $M_\rho$ is also disjoint from all the preimages
$f^{-n}(\pi({\cal V}(F))),\ n\in\N$. This will become important in the proof
of Theorem~\ref{t.sdsm}~.

\end{bem}

%%%%%%%%%%%%%%%%%%%%%%%%%%%%%%%%%%%%%%%%%%%%%%%%%%%%%%%%%%%%%%%%%%
%%%%%%%%%%%%%%%%%%%%%%%%%%%%%%%%%%%%%%%%%%%%%%%%%%%%%%%%%%%%%%%%%%
%%%%%%%%%%%%%%%%%%%%%%%%%%%%%%%%%%%%%%%%%%%%%%%%%%%%%%%%%%%%%%%%%%

\section{Topological Entropy: Proof of Theorem~\ref{t.entropy}} \label{Entropy}

\begin{figure}
\centering
\includegraphics[width=12cm]{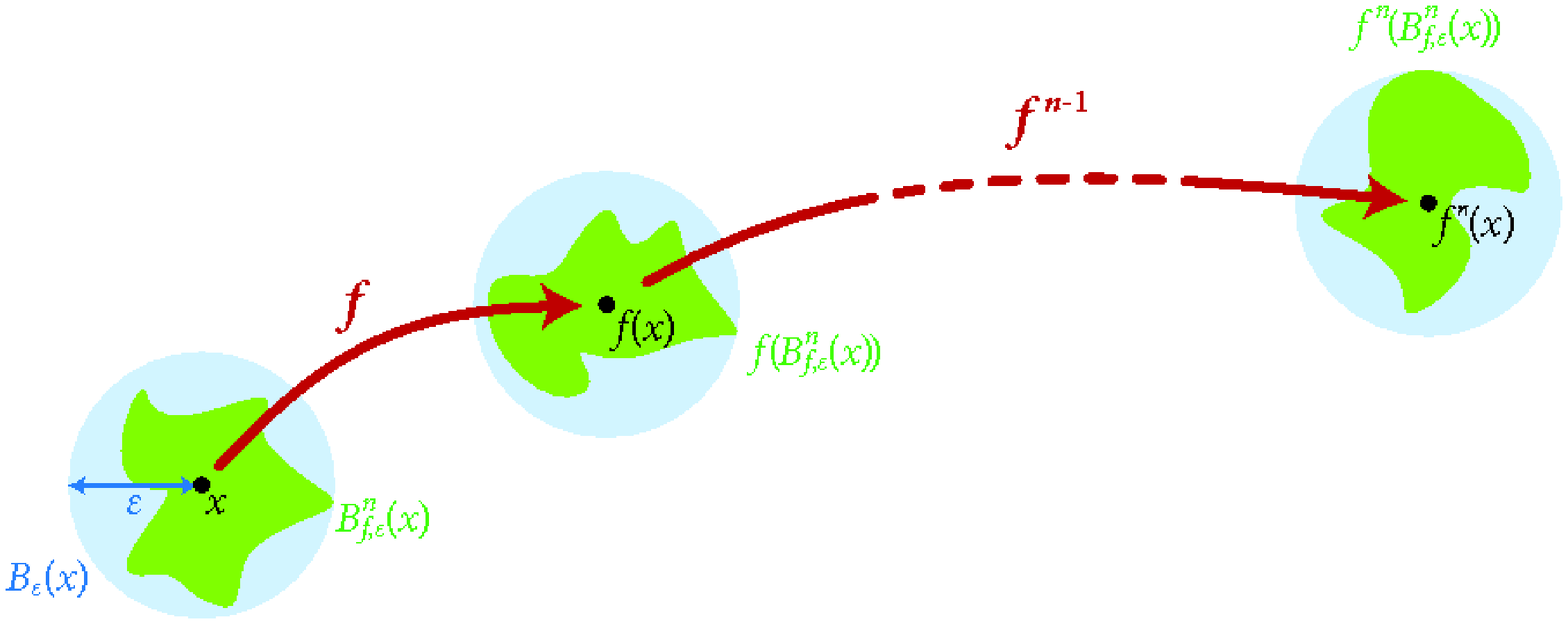}
\caption{Definition of the $(f,n,\eps)$-balls $B^f_{n,\eps}(x)$. }
\label{Def_Eps_Ball_Fig}
\end{figure}

First, we briefly review the definition of topological entropy, following
\cite{bowen:1971} (see also \cite{katok/hasselblatt:1997}). Suppose $(X,d)$ is a
compact metric space and $f:X\to X$ is a continuous map. Then a sequence of
metrics on $X$ is given by $d^f_n(y,z) := \max_{i=0}^n
d(f^i(y),f^i(z))$. For a given $\eps > 0$, $\eps$-balls with respect to $d^f_n$ are called
$(f,n,\eps)${\em -balls} and denoted by $B^f_{n,\eps}(x)$. We let
\[
R(f,n,\eps) \ := \ \min\left\{ k\in\N \mid \exists y_1 \ld y_k \in X : X \ssq
\ikcup B^f_{n,\eps}(y_i)\right\}
\]
and
\[
S(f,n,\eps) \ := \ \max\left\{ k\in\N \mid \exists y_1 \ld y_k \in X :
d^f_n(y_i,y_j) \geq \eps \ \mathrm{ if }\ i \neq j\right\} \ .
\]
We say a finite set $S \ssq X$ is $(f,n,\eps)${\em -separated} if
$d^f_n(y,z) \geq \eps \ \forall y,z\in S: y \neq z$. Then $S(f,n,\eps)$ is the
maximal cardinality of a $(f,n,\eps)$-separated set in $X$. Similarly,
$R(f,n,\eps)$ is the minimal cardinality of a cover of $X$ by
$(f,n,\eps)$-balls. It is easy to see that these quantities are
non-increasing in $\eps$ and satisfy $S(f,n,2\eps) \leq R(f,n,\eps) \leq
S(f,n,\eps)$. Next, we define
$$h_\eps(f) := \limsup_{n\to\infty} \ntel \log R(f,n,\eps)$$ and $$\tilde
h_\eps(f) := \limsup_{n\to\infty} \ntel \log S(f,n,\eps)\ .$$
Again, these two quantities are non-increasing in $\eps$, and the inequalities
$\tilde h_{2\eps}(f) \leq h_\eps(f) \leq \tilde h_\eps(f)$ hold. The topological
entropy of $f$ is defined as
\[
\htop(f) \ := \ \epslim h_\eps(f) \ = \ \sup_{\eps > 0} h_\eps(f) \ ,
\]
and from the preceding discussion it follows that we also have
\[
\htop(f) \ = \ \epslim \tilde h_\eps(f) \ = \ \sup_{\eps > 0} \tilde h_\eps(f)
\ .
\]
We remark that replacing the metric $d$ by another metric $d'$ which is
equivalent (meaning that $d$ and $d'$ induce the same topology) does not
change the topological entropy. In particular, there is no need to specify
below which metric we choose on the product space $\Theta \times \kreis$, any
metric compatible with the product topology will do. However, for simplicity
we will assume the metric on $\Theta \times \kreis$ is chosen such that
$d((\theta_1,x_1),(\theta_2,x_2)) \geq \max\{d(\theta_1,\theta_2),d(x_1,x_2)\}$.

\medskip

For the proof of Theorem~\ref{t.entropy}, it will be convenient to work with a
lift of the original map to a finite covering space of $\Theta \times
\kreis$. Hence, we would like to know that this does not alter the topological
entropy. We denote the $k$-fold cover of the circle by $\kreis_k = \R / k\Z$ and
write $\hat \pi_k$ for the covering map $\hat \pi_k:\Theta \times \kreis_k \to
\Theta \times \kreis$.

\begin{figure}
\centering
\includegraphics[width=8cm]{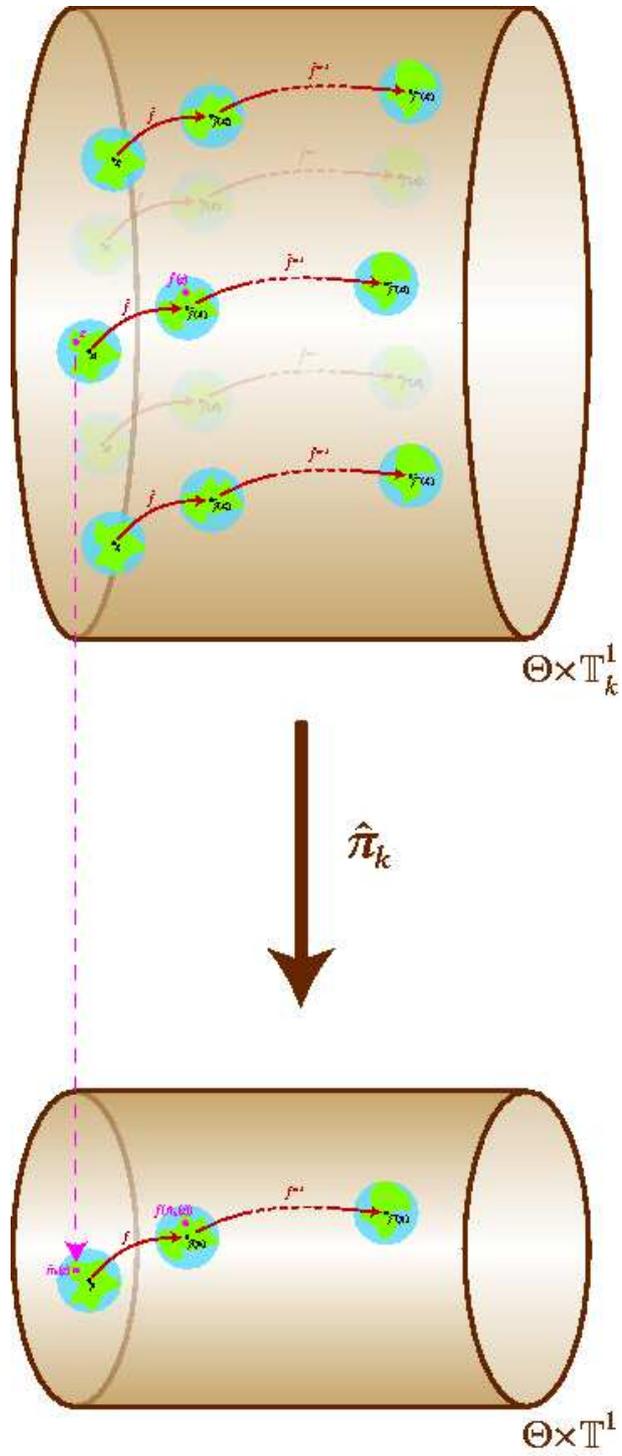}
\caption{Illustration of the proof of Lemma \ref {l.entropy}, with $z \in z_1^2$ }
\label{Covering_Entropy_Fig}
\end{figure}

\begin{lem} \label{l.entropy} Suppose $f:\Theta \times \kreis \to \Theta \times
  \kreis$ is continuous and homotopic to the identity, let $X = \Theta \times
  \kreis_k$ and assume $\hat f : X \to X$ is a lift of $f$. Then $\htop(f) =
  \htop(\hat f)$.
\end{lem}

\proof\ A covering of $X$ with $(\hat f,n,\eps)$-balls projects to a
covering of $\Theta \times \kreis$ with $(f,n,\eps)$-balls. Hence
$R(f,n,\eps) \leq R(\hat f,n,\eps)$ and thus $\htop(f) \leq \htop(\hat
f)$.

In order to prove the converse inequality, let $\eps_0 \in (0,\viertel)$ be such
that $d(y,z) < \eps_0$ implies $d(\hat f(y),\hat f(z)) < \halb$ for all $y,z\in
X$. We will show that for any $0 < \eps \leq \eps_0$ we have $R(\hat f,n,\eps)
\leq k R(f,n,\eps)$, which immediately implies $\htop(\hat f) \leq
\htop(f)$. Fix $0 < \eps \leq \eps_0$, let $R := R(f,n,\eps)$ and choose $y_1
\ld y_R \in \Theta \times \kreis$ such that $\Theta \times \kreis \ssq
\bigcup_{i=1}^R B^f_{n,\eps}(y_i)$. For any $i \in \{1 \ld R \}$, the point
$y_i$ has exactly $k$ lifts $z^j_i$ (Figure \ref{Covering_Entropy_Fig}), with
$d(z^j_i,z^l_i) \geq 1$ whenever $j\neq l$. We claim that $X \ssq
\bigcup_{i=1}^R \bigcup_{j=1}^{k} B^{\hat f}_{n,\eps}(z^j_i)$, so $R(\hat
f,n,\eps) \leq kR$ as required. In order to see this, note that for any $z\in X$
we must have $\hat \pi_k(z) \in B^f_{n,\eps}(y_i)$ for some $i\in\{1\ld R\}$. In
particular $\hat \pi_k(z) \in B_\eps(y_i)$, and therefore $z \in B_\eps(z^j_i)$
for some $j\in\{1\ld k\}$ (Figure \ref{Covering_Entropy_Fig}). Now $\hat
\pi_k(z) \in B^f_{n,\eps}(y_i)$ implies that $\hat f(z)$ is contained in one of
the $k$ $\eps$-balls that make up $(\hat \pi_k)^{-1}B_\eps(f(y_i))$. All of
these are pairwise disjoint and have distance $\geq \halb$ to each other, since
$\eps \leq \eps_0 < \viertel$. Due to the choice of $\eps_0$, we must have $\hat
f(z) \in B_\eps(\hat f(z^j_i))$. By induction on $m$, we thus obtain $\hat
f^m(z) \in B_\eps(\hat f^m(z^j_i))$ for all $m=0 \ld n$. Hence $z \in B^{\hat
  f}_{n,\eps}(z^j_i)$. As $z \in X$ was arbitrary, this completes the proof.

\qed

\proof[Proof of Theorem~\ref{t.entropy}~] Suppose $f$ is a UEF circle
endomorphism, homotopic to the identity. Further, assume $F$ is a lift of $f$
and the rotation interval $\rhofib(F)$ has non-empty
interior. We will work with a lift $\hat f : X \to X$ to the finite covering
space $X= \Theta \times \kreis_4$ and show that the numbers $S(\hat f,n,1)$
grow exponentially.

For any $\rho\in\rhofib(F)$, let $M_\rho$ be the minimal set provided by
Theorem~\ref{t.rotation-interval}.  Choose $\rho_1,\rho_2 \in \rhofib(F)$ with
$\rho_2 > \rho_1$ and let $\epsilon = \frac{1}{4}(\rho_2 - \rho_1)$. Recall that
$ \pi_1(M_{\rho_1}) \cap \pi_1(M_{\rho_2}) \neq \emptyset$
(Remark~\ref{r.projections}). By the uniform convergence of the quantities
$\frac{1}{n}(F^n_\theta(x)-x)$ on $M_{\rho_1}$ and $M_{\rho_2}$ there exists $N
\in \N$ such that for any $\theta \in \pi_1(M_{\rho_1}) \cap \pi_1(M_{\rho_2})$
we have for all $n \geq N$
\begin{align*}
\left\vert F^n_\theta(x_1) - x_1 - n\rho_1 \right\vert &< n\epsilon \\
\left\vert F^n_\theta(x_2) - x_2 - n\rho_2 \right\vert  &< n\epsilon
\end{align*}
for any $x_1,x_2$ such that $(\theta,x_1) \in \pi^{-1}(M_{\rho_1})$ and
$(\theta,x_2) \in \pi^{-1}(M_{\rho_2})$. Thus
\begin{align*}
F^n_\theta(x_1) - x_1  &> n\rho_1 - n\epsilon \\
F^n_\theta(x_2) - x_2  &<  n\rho_2 + n\epsilon
\end{align*}
so that 
$$
F^n_\theta(x_1) - x_1 +  n\rho_2 + n\epsilon > F^n_\theta(x_2) - x_2 + n\rho_1 - n\epsilon
$$
and hence (recall that $\rho_2 - \rho_1 = 4\epsilon$):
\begin{align*}
F^n_\theta(x_1) -  F^n_\theta(x_2)  &> x_1 - x_2 + n(\rho_1 - \rho_2 - 2\epsilon) \\
&> x_1 - x_2 + 2n\epsilon
\end{align*}

\begin{figure}
\centering
\includegraphics[width=12cm]{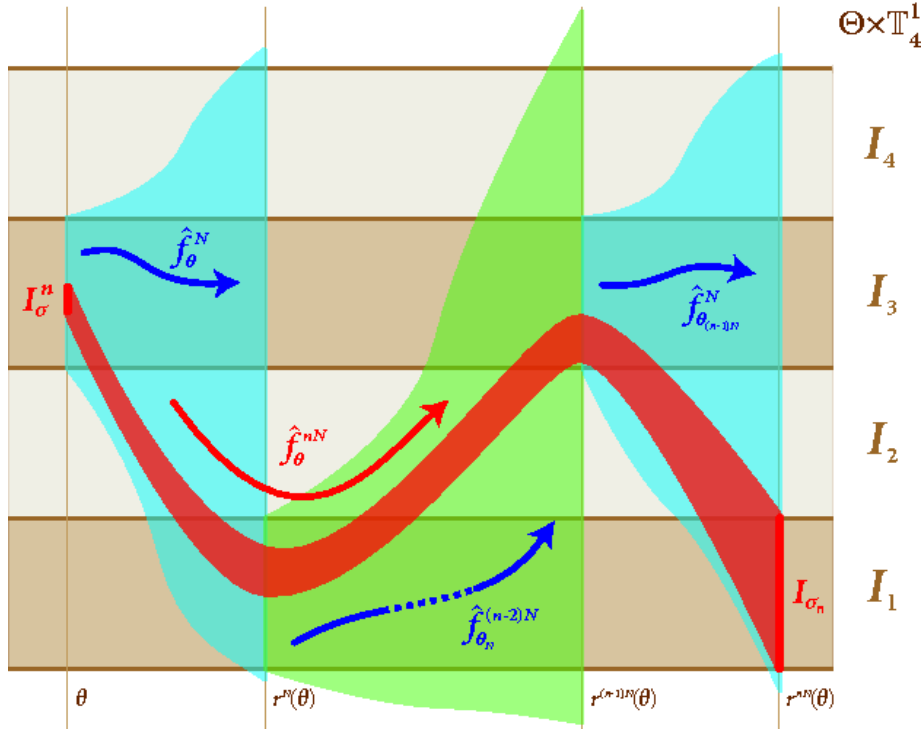}
\caption{Construction of the Set $I^n_\sigma$ defined by
  (\ref{def_partition}). The map $\hat f_\theta^N$ maps any of the intervals
  $I_1 \ld I_4$ at least once around the whole of $\kreis_4$.}
\label{Positive_Entropy_Partition_Fig}
\end{figure}

By \ref{e.periodicity} we can take $x_2 - 1 < x_1 < x_2$ without loss of
generality. We then choose $N$ sufficiently large such that $2N\epsilon > 5$ and
hence
$$
F^N_\theta(x_1) -  F^N_\theta(x_2)  > 4
$$
\begin{figure}
\centering
\includegraphics[width=12cm]{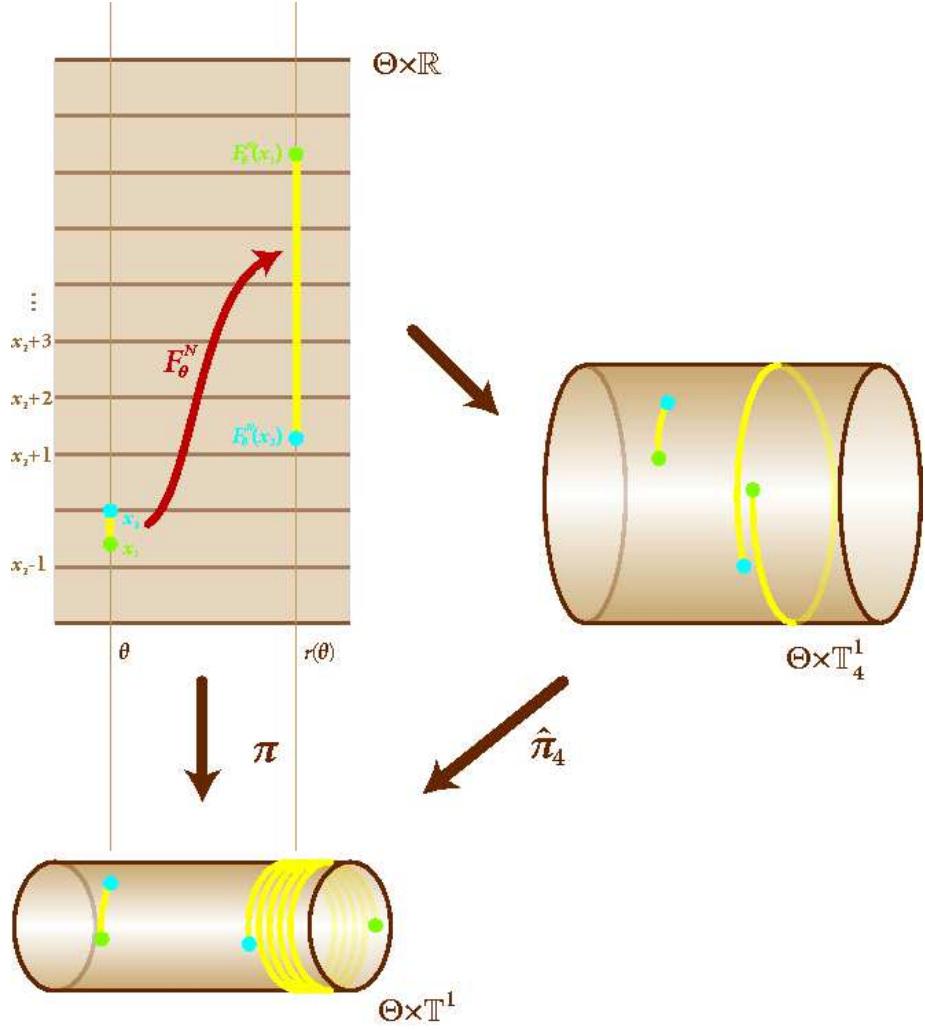}
\caption{Illustration of the proof of Theorem~\ref{t.entropy} }
\label{Positive_Entropy_Proof_Fig}
\end{figure}
for any $x_1,x_2$ such that $(\theta,x_1) \in \pi^{-1}(M_{\rho_1})$ and
$(\theta,x_2) \in \pi^{-1}(M_{\rho_2})$ and $x_2 - 1 < x_1 < x_2$. This implies
that for any $\theta \in \pi_1(M_{\rho_1}) \cap \pi_1(M_{\rho_2})$ the map $\hat
f_\theta^N$ sends each of the intervals $I_i := [i-1,i] \ssq \kreis_4,\ i=1 \ld
4$ surjectively onto $\kreis_4$ (Figure \ref{Positive_Entropy_Proof_Fig}). For
any such $\theta$ and any finite sequence $\sigma \in \{1,3\}^{n+1}, n \in \N$
of the symbols $1$ and $3$ define the set (Figure
\ref{Positive_Entropy_Partition_Fig})
\begin{equation} \label{def_partition}
I^n_\sigma := \cap_{i=0}^{n} (\hat f_\theta^{iN})^{-1}(I_{\sigma_i})
\end{equation}
By definition, $I^0_\sigma$ = $I_{\sigma_0}$ and by the above $I_{\sigma_1} \subset
\hat f_\theta^N(I_{\sigma_0})$. Hence $\hat f^{N}_\theta(I^1_\sigma) \ = \
I_{\sigma_1}$. Similarly $I_{\sigma_2} \subset \hat
f_\theta^N(I_{\sigma_1}) = \hat f_\theta^N(\hat f^{N}_\theta(I^1_\sigma) ) =
\hat f^{2N}_\theta(I^1_\sigma)$ and so $\hat f^{2N}_\theta(I^2_\sigma) =
I_{\sigma_2} $. Continuing by induction we see that
\[
\hat f^{iN}_\theta(I^n_\sigma) \ = \ I_{\sigma_i} \ .
\]
and in particular, $I^n_\sigma$ is non-empty for any $n \in \N$. Clearly, for
any $x \in I^n_\sigma$ and $x'\in I^n_{\sigma'}$ with $\sigma \neq \sigma'$, the
points $(\theta,x)$ and $(\theta,x')$ are $(\hat f^N,n,1)$-separated. Thus
$S(\hat f^N,n,1) \geq 2^{n+1}$. But $S(\hat f,nN,1) \geq S(\hat f^N,n,1)$ and
therefore
\begin{align*}
\tilde h_1(\hat f)  &=  \limsup_{n\to\infty} \ntel \log S(\hat f,n,1)  \\
 & \geq  \limsup_{n\to\infty} \frac{1}{nN} \log 2^{n+1} \ =
\ \frac{\log 2}{N}  > \ 0 \ .
\end{align*}
Since by definition $\htop(\hat f) \geq = h_1(\hat f)$ and $\htop(f) =
\htop(\hat f)$ by Lemma~\ref{l.entropy}, this completes the proof.

\qed

The above proof shows that the positive entropy of $f$ is even realised on
single fibres, meaning that for suitable $\theta\in\Theta$ we can find an
exponentially growing number of $(f,n,\eps)$-separated points contained in
$\{\theta\} \times \kreis$. However, this is by no means surprising. In fact,
when $\htop(r) = 0$, as in the quasiperiodically forced case, it is the only
way to obtain positive topological entropy for the skew-product
transformation. This follows from a well-known result by Bowen. In order to
state it, we have to introduce the topological entropy of a subset $K\ssq X$,
where as at the beginning of this section we assume that $X$ is a compact
metric space. We let
\[
R(f,K,n,\eps) \ := \ \min\left\{ k\in\N \mid \exists y_1 \ld y_k \in X : K
\ssq \ikcup B^f_{n,\eps}(y_i)\right\}
\]
and then define $h_\eps(f,K) := \nLim \ntel \log R(f,K,n,\eps)$ and
$\htop(f,K) := \epslim h_\eps(f,K)$. The numbers $S(f,K,n,\eps)$ and $\tilde
h_\eps(f,K)$ are defined similarly, as above.
\begin{thm}[Bowen \cite{bowen:1971}] \label{t.bowen} Suppose $X,Z$ are compact
  metric spaces, $r:X\to X$, $f: Z\to Z$ and $p:Z \to X$ are continuous maps,
  with $p$ surjective and $p\circ f = r \circ p$. Then \[ \htop(f) \ \leq \
  \htop(r) + \sup_{y\in X} \htop(f,p^{-1}\{y\}) \ . \]
\end{thm}
Hence, if $f$ is a UEF circle endomorphism and $\htop(r)=0$, then $\htop(f) >
0$ implies that there exists some $\theta \in \Theta$ with $\htop(f,\{\theta\}
\times \kreis) > 0$. Conversely, if all fibre maps are monotone then the above
theorem easily entails the following
\begin{cor} \label{c.bowen} Suppose $f$ is a UEF monotone circle map. Then
  $\htop(f) = \htop(r)$. In particular, if $f$ is a QPF monotone circle map,
  then $\htop(f) = 0$. Note that $f$ need not necessarily be homotopic to the
  identity.
\end{cor}
\proof For any $\theta \in \Theta$, let
$\T_\theta := \{\theta\} \times \kreis$. In view of Theorem~\ref{t.bowen},
we only have to prove that
\[
\htop(f,\T_\theta) = 0 \quad \forall \theta \in \Theta \ .
\]
In order to do so, we will show that the numbers
$R(f,\T_\theta,n,\eps)$ can grow at most linearly with $n$. To
that end, fix $\theta_0 \in \Theta$ and $\eps > n$. By compactness, there
exists a finite cover of $\Theta \times \kreis$ by boxes $A^i_j = A_i \times
I_j$, $(i,j=1\ld N)$, with the following properties: {\em \romanlist
\item each set $A^i_j$ has diameter less then $\eps$;
\item there exist $a_0 < a_1 < \ld < a_N=a_0 \in \kreis$, such that $I_j = [a_{j-1},a_j]$;
\item for all $m\in\N$ and $j=1\ld N$, the point $a_j$ has a unique preimage
  under the map $f^m_{\theta_0}$; \listend }
Concerning (iii) note that, due to monotonicity, for each fibre map $f_\theta$
the set of points on which $f_\theta$ is not injective is an at most countable
union of intervals, and each of these intervals is mapped to a single
point. Consequently, for any $m$ there is an at most countable set $E_m$ of
exceptional points, whose preimage under $f^m_{\theta_0}$ is not unique. It
suffices to choose the $a_j$ outside the resulting countable union $\mcup
E_m$.

We denote by ${\cal A}^n$ the $n$-th refinement of the cover ${\cal A}=\{A^i_j
\mid i,j \in \{1\ld N\}\}$, that is
\[
{\cal A}^n \ := \ \left\{ \alpha \ssq \Theta \times \kreis \left| \ \alpha =
\bigcap_{k=0}^n f^{-k}\left(A^{i_k}_{j_k}\right),\ i_k,j_k \in \{1\ld N\} \
\forall k=0\ld n \right. \right\} \ .
\]
By ${\cal A}^n_{\theta_0}$, we denote the restriction of ${\cal A}^n$ to the
$\theta_0$-fibre, that is
\[
{\cal A}^n_{\theta_0} \ := \ \left\{ \alpha \cap \T_{\theta_0} \mid
\alpha \in {\cal A}^n \right\} \smin \{\emptyset\} \ .
\]
Choose points $x_\beta \in \Theta\times \kreis$, such that $x_\beta \in \beta
\ \forall \beta \in {\cal A}^n_{\theta_0}$. Since the sets $A^i_j$ all have
diameter less than $\eps$, we obtain
\[ \T_{\theta_0} \ \ssq \ \bigcup_{\beta\in {\cal A}^n_{\theta_0}} B^f_{n,\eps}(x_\beta) \ .
\]
Consequently, $R(f,\T_{\theta_0},n,\eps) \leq \# {\cal A}^n_{\theta_0}$ where
$\# {\cal A}^n_{\theta_0}$ denotes the number of elements in the partition $
{\cal A}^n_{\theta_0}$ . However, if $k\in\{0\ld n\}$ is fixed, then due to
montonicity and properties (ii) and (iii) above the preimages of the intervals
$I_j$ under the map $f^k_{\theta_0}$ are all intervals with pairwise disjoint
interiors. It follows that ${\cal A}^n_\theta$ is a partition of the circle
$\T_{\theta_0}$, given by the points $(\theta_0,(f^k_{\theta_0})^{-1}(a_j))$,
$j=1 \ld N, k=0 \ld n$. This implies that $\# {\cal A}^n_\theta \leq (n+1)N$,
and consequently $R(f,\T_{\theta_0},n,\eps) \leq (n+1)N$. Since $\theta_0 \in
\Theta$ was arbitrary, this completes the proof.

\qed

\proof[\bf Proof of Theorem~\ref{t.rotation-interval} (ii) ] This follows
immediately from Corollary~\ref{c.bowen}. Recall from the proof of part (i) that
the orbit under $F$ of any point in $M_\rho$ coincides with the orbit under the
UEF monotone map $F_{t(\rho)}$. Hence $\htop(f|_{M_\rho}) \leq
\htop(f_{t(\rho)})$, where $F_{t(\rho)}$ is a lift of $f_{t(\rho)}$. But, since
$F_{t(\rho)}$ is monotone we have $\htop(f_{t(\rho)}) = 0 $, as required.

\qed

%%%%%%%%%%%%%%%%%%%%%%%%%%%%%%%%%%%%%%%%%%%%%%%%%%%%%%%%%%%%%%%%%%
%%%%%%%%%%%%%%%%%%%%%%%%%%%%%%%%%%%%%%%%%%%%%%%%%%%%%%%%%%%%%%%%%%
%%%%%%%%%%%%%%%%%%%%%%%%%%%%%%%%%%%%%%%%%%%%%%%%%%%%%%%%%%%%%%%%%%

\section{Strangely Dispersed Minimal Sets: Proof of Theorem~\ref{t.sdsm}} \label{SDMS}

\begin{figure}
\centering
\includegraphics[width = 12cm]{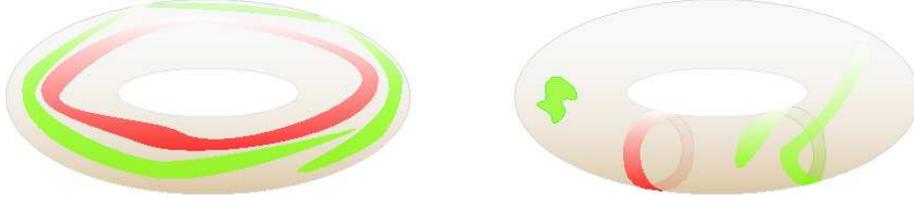}
\caption{Essentially bounded sets on the torus. Green sets are essentially bounded, red sets are not.}
\label{Essentially_Bounded_Fig}
\end{figure}

Let $q:\R^2 \to \ntorus = \R^2 / \Z^2$ denote the quotient map. We call a subset
$E \ssq \ntorus$ {\em essentially bounded}, if all connected components of
$q^{-1}(E)$ are bounded, Figure \ref{Essentially_Bounded_Fig}. The following
proposition will be the key ingredient in the proof of Theorem~\ref{t.sdsm}~.

\begin{prop} \label{p.sdms-criterion}
Suppose $f$ is a QPF circle endomorphism, homotopic to the identity. Further,
assume $E\ssq \ntorus$ is open and essentially bounded and $M \ssq E$ is a
minimal set. Then $M$ is strangely dispersed.
\end{prop}
\proof\ Since $M$ is minimal by assumption, it remains to show that it has
properties (ii) and (iii) in Definition~\ref{d.sdms}~.

In order to see that connected components of $M$ are contained in single fibres,
suppose that $\hat E_0$ is a connected component of $\hat E := q^{-1}(E)$ such
that $\hat M_0 := q^{-1}(M) \cap \hat E_0 \neq \emptyset$. Since $E$ is
essentially bounded, $\hat E_0$ is bounded and hence $\hat M_0$ is compact. Thus
the first coordinate of points in $\hat M_0$ attains a minimal value, at say
$(\hat\theta_0,\hat x_0) \in \hat M_0$ (more precisely $\hat \theta_0 =
\inf\{\theta\in \R \mid \exists x \in \R : (\theta,x) \in \hat M_0\})$. Let
$(\theta_0,x_0) := q(\hat \theta_0,\hat x_0)$ and $E_0 := q(\hat E_0)$. Observe
that $E_0\ssq E$ is a connected component of $E$ and in particular, $E_0$ is
open. Also $(\theta_0,x_0) \in M$ and $(\theta_0,x_0) \in E_0$

Now, assume that $C \ssq M$ is a connected component of $M$ which is not
contained in a single fibre. Then $\pi_1(C)$ is connected and hence an interval
of positive length, say $\pi_1(C) = [a,b]$ with $\delta = d(a,b)> 0$. We assume
for simplicity of exposition that $d(a,b) < \halb$. Choose $(\theta,x) \in C$
with $d(\theta,a) = d(\theta,b) = \delta/2$. Observe that for any $n\in\N$, we
have $\pi_1(f^n(C)) = [r^n(a),r^n(b)]$ which also has length $\delta$, and
$d(\theta_n,a) = d(\theta_n,b) = \delta/2$ where $(\theta_n,x_n) =
f^n(\theta,x)$.

Since $M$ is minimal, the orbit of $(\theta,x)$ is dense in $M$. By the above,
$E_0$ is open and contains $(\theta_0,x_0) \in M$, so that there exists some
$n\in\N$, such that $f^n(\theta,x) \in E_0 \cap B_{\delta/4}(\theta_0,x_0)$. The
set $f^n(C)$ is connected and $f^n(C) \ssq M \ssq E$ for all $n\in\N$. Hence
$f^n(C)$ is contained in a connected component of $E$. Since $f^n(C)$ contains
$f^n(\theta,x) \in E_0$ this connected component must by $E_0$, that is $f^n(C)
\ssq E_0$.

Define $\hat D_0$ as the unique connected component of $q^{-1}(f^n(C))$ that
contains the unique point $(\hat \theta^*,\hat x^*)$ in $q^{-1}\{f^n(\theta,x)\}
\cap B_{\delta/4}(\hat \theta_0,\hat x_0)$. Then $(\hat \theta^*,\hat x^*) \in
\hat E_0$, and by connectedness $\hat D_0 \ssq \hat E_0$ and hence $\hat D_0
\ssq \hat M_0$. Since $q(\hat D_0) = f^n(C)$, $q(\hat \theta^*,\hat x^*) =
f^n(\theta,x) = (\theta_n,x_n)$ and $\pi_1(f^n(C)) = [\theta_n - \delta/2,
\theta_n + \delta/2]$ we have $\pi_1(\hat D_0) = [\hat \theta^* - \delta/2, \hat
\theta^* + \delta/2]$.

But recall that $\hat \theta_0 = \inf\{\theta\in \R \mid \exists x \in \R :
(\theta,x) \in \hat M_0\})$ and since $\hat D_0 \ssq \hat M_0$ we must have
$\hat \theta_0 \ \leq \ \hat \theta^* - \delta/2$. On the other hand $(\hat
\theta^*,\hat x^*) \in B_{\delta/4}(\hat \theta_0,\hat x_0)$, so that $\hat
\theta^* \leq \hat \theta_0 + \delta/4$ which implies that $\hat \theta^* -
\delta/2 \leq \hat \theta_0 - \delta/4$. Combining these two inequalities yields
the contradiction
\[
\hat \theta_0 \ \leq \ \hat \theta^* - \delta/2 \ \leq \ \hat \theta_0 - \delta/4 \ .
\]
Hence any connected component of $M$ must be contained in a single fibre,
proving property (ii).

It remains to prove that for any $\thx \in M$ and any open neighbourhood $U$ of
$\thx$, the set $\pi_1(U\cap M)$ contains a non-empty open interval,
i.e. property (iii). First observe that if this property holds for some $\thx
\in M$ then it holds for $f\thx$ (and hence $f^n\thx$ for any $n \in \N$). To
see this, let $U$ be an open neighbourhood of $\thx$. Then $f^{-1}(U)$ is an
open neighbourhood of $\thx$, and hence $\pi_1(f^{-1}(U) \cap M)$ contains a
non-empty open interval $(a,b)$. Hence $\pi_1(U\cap M)$ contains $(r(a),r(b))$
which has the same length as $(a,b)$ and hence is a non-empty open
interval. Also, property (iii) is closed, that is if it holds for a convergent
sequence of points $(\theta_i,x_i) \in M$ with $(\theta_i,x_i) \to \thx$ then it
holds for the limit point $\thx$. This is because if $U$ is an open
neighbourhood $U$ of $\thx$ then it is an open neighbourhood of $(\theta_i,x_i)$
for all sufficiently large $i$.

Thus, property (iii) is both closed and invariant and hence, it either holds for
all or for no point in $M$ since by minimality the only closed invariant subsets
of $M$ are the empty set and $M$ itself. Arguing by contradiction, let us assume
that every $z \in M$ has a neighbourhood $U(z)$, such that $\pi_1(U(z) \cap M)$
contains no open interval and hence is nowhere dense. By compactness, $M$ is
covered by a finite number $U(z_1) \ld U(z_N)$ of such neighbourhoods. However,
this would imply that $\pi_1(M)$ is the union of a finite number of nowhere
dense sets and hence is itself nowhere dense. This is clearly a contradiction,
since $\pi_1(M)$ must be the whole circle, because this is the only closed
invariant set of the underlying irrational rotation.

\qed

\proof[Proof of Theorem~\ref{t.sdsm}] Suppose $f$ is given by
(\ref{e.arnold}), and consequently has a lift $F$ with fibre maps
\begin{equation*} \label{e.arnold-lift}
F_\theta(x) \ = \ x + \tau + \frac{\alpha}{2\pi} \sin(2\pi x) + \beta
\sin(2\pi\theta) \ .
\end{equation*}

\begin{figure}
\centering
\subfloat[The set ${\cal V}(F_t)$]{
\includegraphics[width = 6cm]{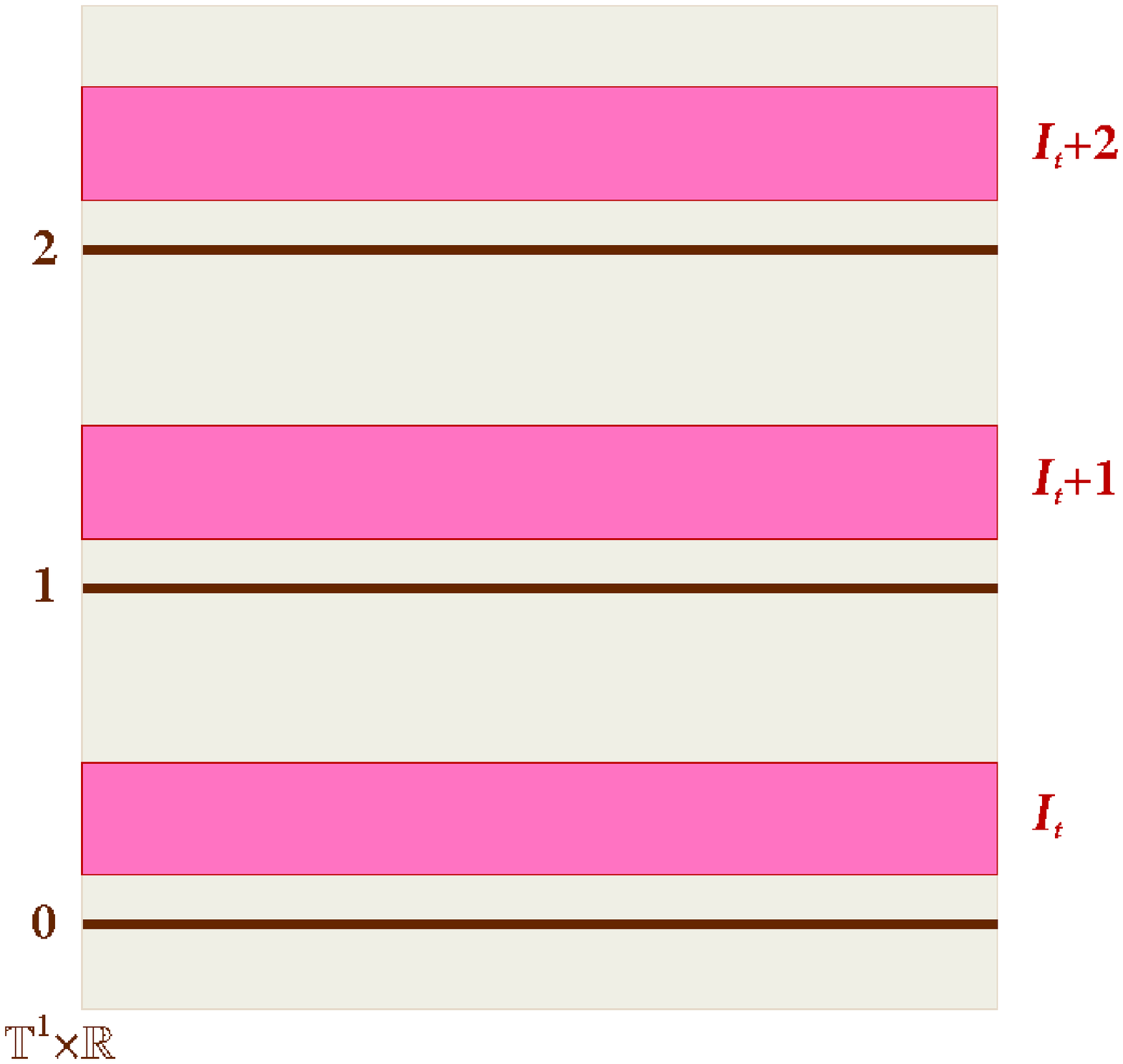}
}
\quad
\subfloat[The sets $ \kreis \times I'$ and $F_t^{-1}(\kreis \times I')$.]{
\includegraphics[width = 6cm]{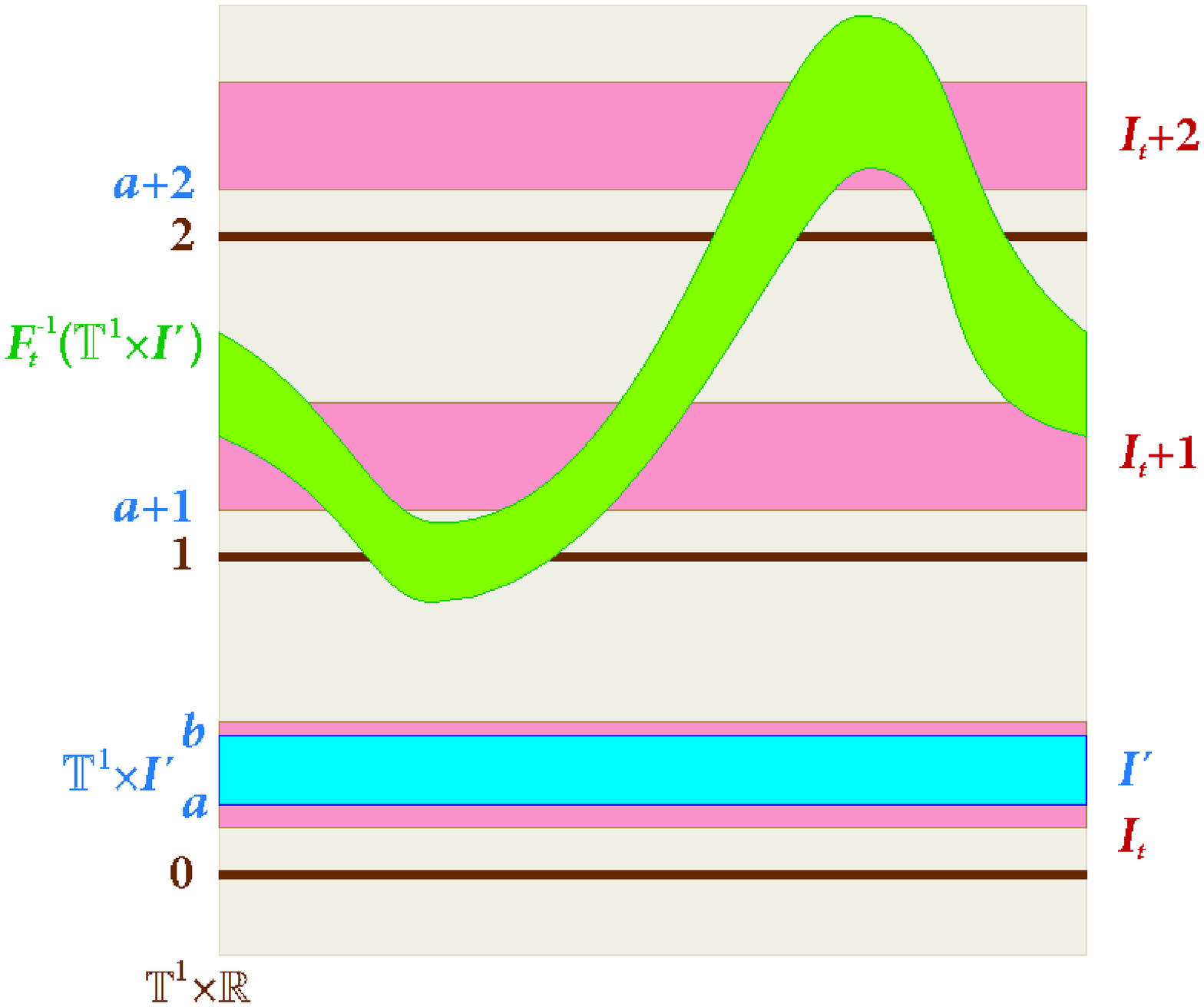}
}\\
\subfloat[The set $\hat E$, coloured yellow, which is the lift of $E = \ntorus
\smin \left( \kreis \times \pi(I') \cup f_t^{-1}(\kreis \times \pi(I'))
\right)$]{
\includegraphics[width = 6cm]{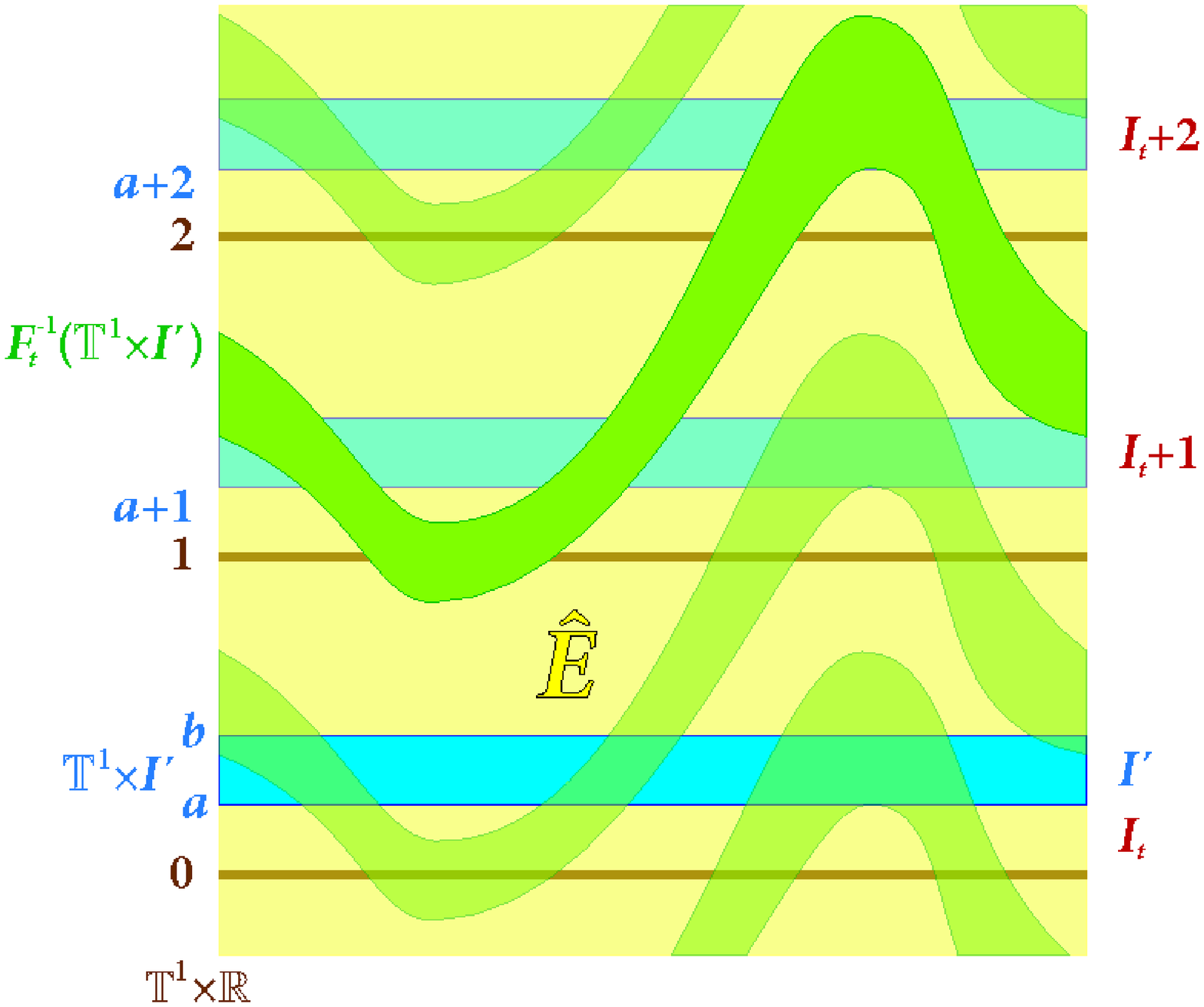}
} \quad \subfloat[The curve $\Gamma$ in $F_t^{-1}(\kreis \times I')$ and its
image under $F_t$ in $ \kreis \times I'$.]{
\includegraphics[width = 6cm]{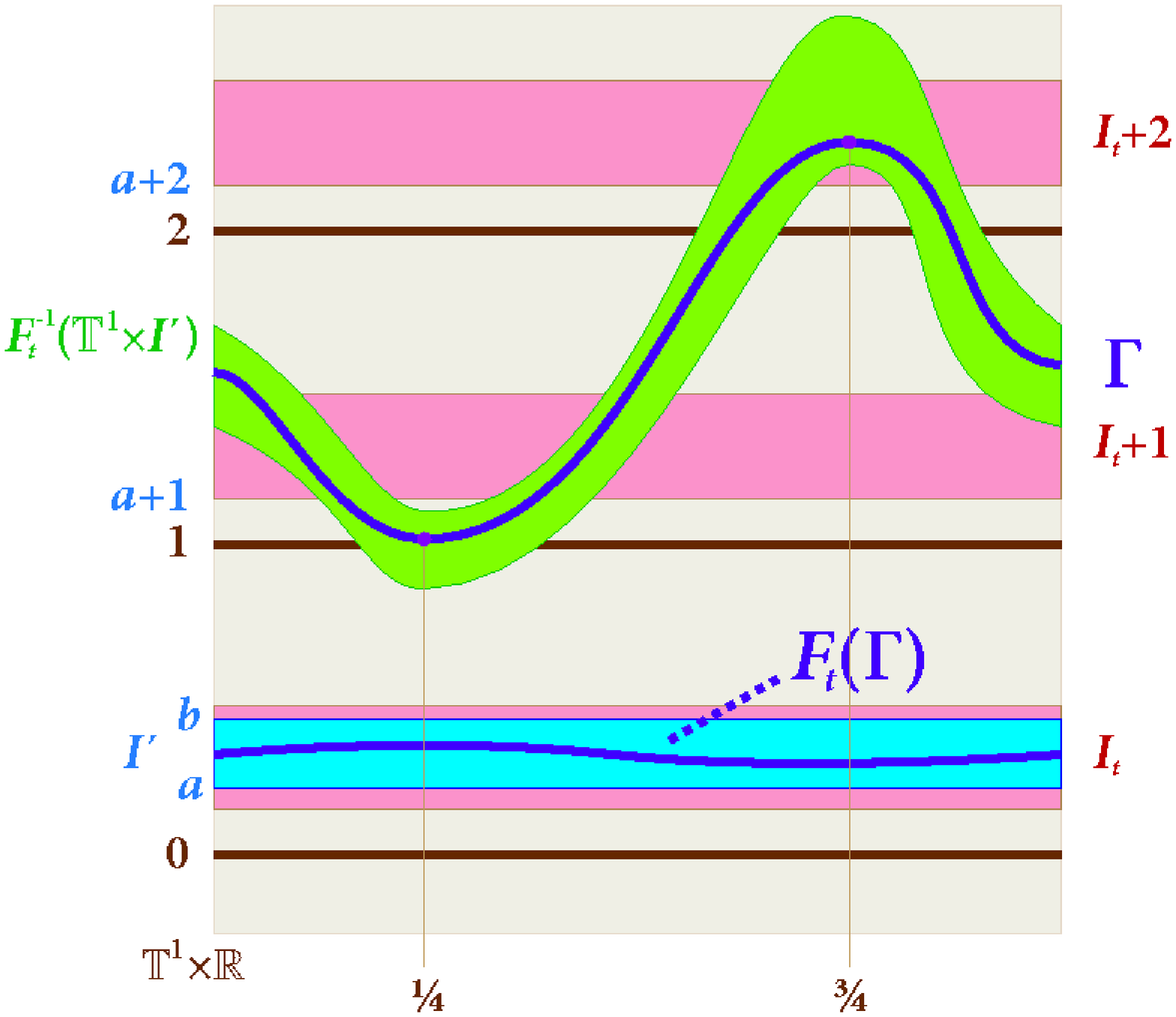}
\label{Plateau_Strips_Gamma_Fig}
}
\caption{Proof of Theorem~\ref{t.sdsm}. Construction of the set $E$. This is
  given by the projection to \ntorus of the complement, coloured yellow, of $
  \kreis \times I'$ (blue) and $F_t^{-1}(\kreis \times I')$ (green) in b).}
\label{Plateau_Strips_Fig}
\end{figure}

{\em Part (a).} Recall that the QPF plateau maps $F_t$ in the proof of
Theorem~\ref{t.rotation-interval} were given by $F_{t,\theta} :=
(F_\theta)_t$, with the mapping $[0,1] \times {\cal E},\ (t,G) \mapsto G_t$
provided by Proposition~\ref{p.plateau-homotopy}~. Any $F_t$ induces a QPF
monotone circle map, which we will denote by $f_t$. If we let
\[
G(x) \ := \ x + \tau + \frac{\alpha}{2\pi} \sin(2\pi x) \ ,
\]
then $F_{t,\theta}(x) = G_t(x) + \beta\sin(2\pi\theta)$. In particular, the
plateaus of $F_{t,\theta}$ do not depend on $\theta$. The fact that $f$ in
(\ref{e.arnold}) is bimodal further implies that these plateaus are unique
modulo addition of integers, that is ${\cal U}(G_t) = \bigcup_{n\in\Z} I_t+n$
for some interval $I_t \ssq \R$. Hence, recalling Remark \ref{r.avoid-plateaus},
we have (Figure \ref{Plateau_Strips_Fig}a))
\[ {\cal V}(F_t) \ = \ \bigcup_{\theta \in \kreis} \{\theta\} \times {\cal
  U}(F_{t,\theta}) \ = \ \bigcup_{k\in\Z} \kreis \times (I_t+k) \ .
\]

Now suppose that $\rho \in \rhofib(F)$ and denote by $M_\rho$ the minimal set
obtained in the proof of Theorem~\ref{t.rotation-interval}. Let $t=t(\rho) \in
[0,1]$ be the corresponding parameter such that $\rho(F_t) = \rho$ and $M_\rho$
is a $F_t$-minimal set. As mentioned in Remark~\ref{r.avoid-plateaus}, $M_\rho$
is disjoint from the set $\bigcup_{n\in\N} f^{-n}(\pi({\cal V}(F_t))$. Now if
$I'=[a,b] \ssq I_t$ is a closed interval let 
\[
E \ := \ \ntorus \smin \left( \kreis \times \pi(I') \cup f_t^{-1}(\kreis \times
  \pi(I')) \right) \ ,
\]
as indicated in Figure \ref{Plateau_Strips_Fig}(b,c).  Then $M_\rho \ssq E$, and
in view of Proposition~\ref{p.sdms-criterion} we only have to show that $E$ is
essentially bounded. Since the complement $q^{-1}(E)^c$ of $q^{-1}(E)$ contains
the horizontal line $\R \times \{a\}$ and all its integer translates, it is
obvious that all connected components of $q^{-1}(E)$ are bounded in the vertical
direction.

We also claim that $q^{-1}(E)^c$ contains a continuous curve joining $\R \times
\{a\}$ and $\R \times \{a+1\}$, which implies immediately that it is bounded
horizontally. We in fact show that $V := \pi^{-1}(E)$ contains a continuous
curve joining $\kreis \times \{a+k\}$ and $\kreis \times \{a+k+1\}$ for some
$k\in \Z$. Observe that $F_t^{-1}(\kreis \times I')$ contains a curve $\Gamma$
that is the graph $\Gamma := \{(\theta,\gamma(\theta)\mid \theta \in \kreis\}$
of a continuous function $\gamma : \kreis \to \R$ (Figure
\ref{Plateau_Strips_Gamma_Fig}). Since $\Gamma$ is mapped into $\kreis \times
I'$ we have
\[
F_{t,\theta}(\gamma(\theta)) \ = \ G_t(\gamma(\theta)) + \beta\sin(2\pi\theta)
\ \in I' \quad \forall \theta\in\kreis \ .
\]

Since we assume that $\beta \geq \frac{3}{2}$ we have
\begin{align*}
F_{t,\theta}(\gamma(1/4)) &=  G_t(\gamma(1/4))) + \beta \geq G_t(\gamma(1/4)) + \frac{3}{2} \\
F_{t,\theta}(\gamma(3/4)) &=  G_t(\gamma(3/4)) - \beta \leq G_t(\gamma(3/4)) - \frac{3}{2}.
\end {align*}
Since the length of $I'$ is less than 1 this implies 
$$
1 >  G_t(\gamma(1/4)) +  \frac{3}{2} -  G_t(\gamma(3/4)) +  \frac{3}{2}
$$
and hence
$$
G_t(\gamma(3/4))  - G_t(\gamma(1/4)) >  2.
$$
Since $G_t$ is monotone and $G_t(x+n) = G_t(x)+n \ \forall n\in\Z, $ if
$\gamma(3/4) \leq \gamma(1/4) + 2$ then $G_t(\gamma(3/4)) \leq G_t(\gamma(1/4) +
2) = G_t(\gamma(1/4)) + 2$, so that $G_t(\gamma(3/4)) - G_t(\gamma(1/4)) \leq
2$. Hence we must have $\gamma(3/4) > \gamma(1/4) + 2$, or in other words
\[
\gamma(1/4) - \gamma(3/4) \geq 2.
\]
Hence there exists $k\in\Z$, such that $\Gamma$ intersects both $\kreis \times
\{a+k\}$ and $\kreis \times \{a+k+1\}$. This proves our claim.

\begin{figure}
\centering
\includegraphics[width = 8cm]{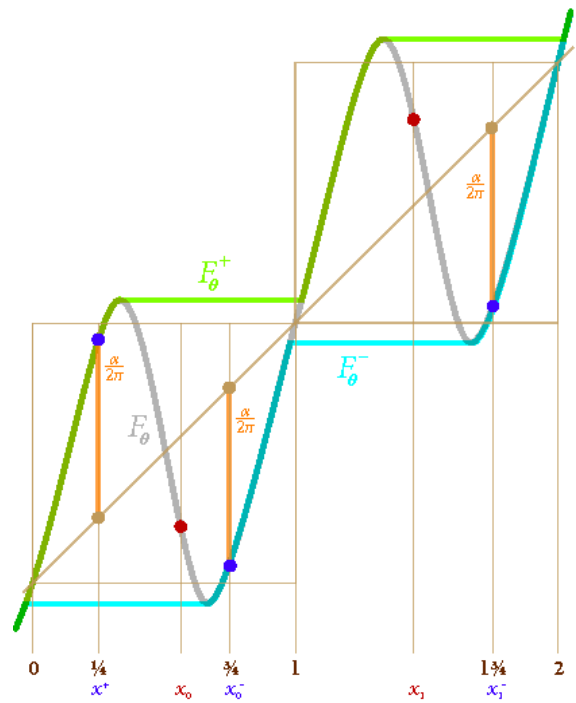}
\caption{Proof of the fact that $F^+_\theta(x) \geq F^-_\theta(x) + \alpha/\pi -
  (x^- - x^+)$ We consider two cases: one where $1/4 \leq x \leq 3/4$ (shown as
  $x_0$) and the other where $3/4 \leq x \leq 7/4$ (shown as $x_1$). For both of
  these $x^+ = 1/4$, whereas $x^- = 3/4$ and $7/4$ respectively, indicated as
  $x^-_0$ and $x^-_1$ on the figure.}
\label{Positive_Rotation_Interval_Fig}
\end{figure}

{\em Part (b)} Recall from the proof of Theorem~\ref{t.rotation-interval} that
$\rhofib(F) = [\rho_1,\rho_2]$, where $\rho_1 = \rho(F^-)$ and $\rho(F_2) =
\rho(F^+)$. For any $x\in\R$, let $x^-:= \inf\{y \in \Z + \dreiviertel \mid y
\geq x\}$ and $x^+ := \sup\{ y \in \Z + \viertel \mid y \leq x\}$ (Figure
\ref{Positive_Rotation_Interval_Fig}). Then $x^- = x^+ + \halb$ or $x^- = x^+ +
\frac{3}{2}$. Note that
\begin{align*}
F_\theta(x^+) &= x^+ + \tau + \frac{\alpha}{2\pi} + \beta \sin(2\pi \theta) \\
F_\theta(x^-)  &= x^- + \tau - \frac{\alpha}{2\pi} + \beta \sin(2\pi \theta)
\end {align*}
and hence
$$
F_\theta(x^+) - F_\theta(x^-) = x^+ - x^- +  \frac{\alpha}{\pi}.
$$

Recall from the defintion of plateau maps that if $x' \leq x$ then
$F^-_\theta(x') \leq F(x)$ and $F^+_\theta(x) \geq F(x')$. Since $x^+ \leq x
\leq x^-$ we have
\begin{align*}
F^+_\theta(x) & \geq F_\theta(x^+) \\
F_\theta(x^-)  &\geq \ F^-_\theta(x)
\end {align*}
for all $\thx \in\kreis \times \R$. Thus if $\alpha \geq \frac{5}{2}\pi$, it
follows that
\[
F^+_\theta(x) \ \geq F_\theta(x^+) \ \geq \ F_\theta(x^-) + \frac{5}{2} - (x^- -
x^+) \ \geq \ F^-_\theta(x) +1 \ .
\]
Hence we have $F^+_\theta(x) \geq F^-_\theta(x) +1$ for all $ \thx \in\kreis
\times \R$, which implies $\rho_2 \geq \rho_1 + 1$. Hence $\rhofib(F) =
[\rho_1,\rho_2]$ has positive length, as required.

\qed

%\bibliography{snaphysics,qpfs,torus,dynamics} \bibliographystyle{unsrt}

\end{document}